\documentclass[
submission
]{dmtcs-episciences}

\usepackage[utf8]{inputenc}
\usepackage{subfigure}
\usepackage{tikz}
\usepackage{amsfonts,amsmath,amssymb}

\newcommand{\pch}{\chi_{\rho}}
\newcommand{\diam}{{\rm diam}}
\def\cp{\,\square\,}

\newtheorem{theorem}{\bf Theorem}[section]
\newtheorem{corollary}[theorem]{\bf Corollary}
\newtheorem{obser}[theorem]{\bf Observation}

\newtheorem{lemma}[theorem]{\bf Lemma}
\newtheorem{proposition}[theorem]{\bf Proposition}

\tikzstyle{vertex}=[circle, draw, inner sep=0pt, minimum size=6pt]

\usepackage[round]{natbib}

\author{Sandi Klav\v zar\affiliationmark{1,2,3}\thanks{The financial support of the Slovenian Research Agency (research core funding No.\ P1-0297,  bilateral grant BI-US/18-19-011, project J1-9109) and of the Wylie Enrichment Fund of Furman University is acknowledged.}
  \and Douglas F. Rall\affiliationmark{4}\thanks{The financial support of the Slovenian Research Agency (bilateral grant BI-US/18-19-011) and of the Wylie Enrichment Fund of Furman University is acknowledged.}
  }
\title[Packing chromatic vertex-critical graphs]{Packing chromatic vertex-critical graphs}
\affiliation{
  Faculty of Mathematics and Physics, University of Ljubljana, Slovenia\\
  Faculty of Natural Sciences and Mathematics, University of Maribor, Slovenia\\
  Institute of Mathematics, Physics and Mechanics, Ljubljana, Slovenia\\
  Department of Mathematics, Furman University, Greenville, SC, USA}
\keywords{packing chromatic number; packing chromatic vertex-critical graph; tree; caterpillar; Cartesian product of graphs, vertex-transitive graph}
\received{2018-10-10}
\revised{2019-1-31}
\accepted{2019-2-7}
\begin{document}
\publicationdetails{21}{2019}{3}{8}{4878}
\maketitle
\begin{abstract}
The packing chromatic number $\pch(G)$ of a graph $G$ is the smallest integer $k$ such that the vertex set of $G$ can be partitioned into sets $V_i$,  $i\in [k]$, where vertices in $V_i$ are pairwise at distance at least $i+1$. Packing chromatic vertex-critical graphs, $\pch$-critical for short, are introduced as the graphs $G$ for which $\pch(G-x) < \pch(G)$ holds for every vertex $x$ of $G$. If $\pch(G) = k$, then $G$ is $k$-$\pch$-critical. It is shown that if $G$ is $\pch$-critical, then the set $\{\pch(G) - \pch(G-x):\ x\in V(G)\}$ can be almost arbitrary. The $3$-$\pch$-critical graphs are characterized, and $4$-$\pch$-critical graphs are characterized in the case when they contain a cycle of length at least $5$ which is not congruent to $0$ modulo $4$. It is shown that for every integer $k\ge 2$ there exists a $k$-$\pch$-critical tree and that a $k$-$\pch$-critical caterpillar exists if and only if $k\le 7$. Cartesian products are also considered and in particular it is proved that if $G$ and $H$ are vertex-transitive graphs and $\diam(G) + \diam(H) \le \pch(G)$, then $G\cp H$ is $\pch$-critical.
\end{abstract}


\section{Introduction}

Given a graph $G$ and a positive integer $i$, an {\em $i$-packing} in $G$ is a subset $W$ of the vertex set of $G$ such that the distance between any two distinct vertices from $W$ is greater than $i$. Note that in this terminology independent sets are precisely $1$-packings. The {\em packing chromatic number} of $G$ is the smallest integer $k$ such that the vertex set of $G$ can be partitioned into sets $V_1,\ldots, V_k$, where $V_i$ is an $i$-packing for each $i\in [k] = \{1,\ldots, k\}$. This invariant is well defined on any graph $G$ and is denoted $\pch(G)$. More generally, for a nondecreasing sequence $S=(s_1,\ldots,s_k)$ of positive integers, the mapping $c:V(G)\longrightarrow [k]$ is an {\em $S$-packing coloring} if for any $i$ in $[k]$ the set $c^{-1}(i)$ is an $s_i$-packing, see~\cite{goddard-2012}.

The packing chromatic number was introduced  in~\cite{goddard-2008} under the name broadcast chromatic number; the present name is used since~\cite{bkr-2007}. Among the many developments on the packing chromatic number we point out  the very exciting development on the class of subcubic graphs.  The question whether in subcubic graphs the packing chromatic number is bounded by an absolute constant was asked in the seminal paper~\cite{goddard-2008}. \cite{gt-2016} were the first to find an explicit cubic graph $G$ with large packing chromatic number, more precisely $\pch(G)=13$, and asked whether $13$ is an actual upper bound for $\pch$ in the class of cubic graphs. That this is not the case was first demonstrated in~\cite{bkrw-2017a}. Then, in a breakthrough result, \cite{balogh-2018} proved that almost every $n$-vertex cubic graph of girth at least $g$ has the packing chromatic number greater than $k$. \cite{bresar-2018a} complemented the result with an explicit infinite family of subcubic graphs with unbounded packing chromatic number. The packing chromatic number of subcubic outerplanar graphs has been thoroughly investigated in~\cite{bresar-2018+} and in \cite{gastineau-2018}, while for the packing chromatic number of some additional classes of cubic graphs see~\cite{laiche-2018}. From the other investigations on the packing chromatic number we emphasize that the decision version of the packing chromatic number is NP-complete on trees~\cite{fiala-2010} and that the packing chromatic number of the infinite square lattice is from the set $\{13, 14, 15\}$~\cite{barnaby-2017}. For further results on the packing chromatic number see the very recent paper~\cite{korze-2019} and references therein. 

An important consideration when investigating a graph invariant is to understand critical graphs with respect to the invariant. The most prominent example is probably the concept of vertex-critical graphs with respect to the chromatic number; that is, graphs $G$ such that $\chi(G-x)<\chi(G)$ for every vertex $x$ of $G$, cf.~\cite{jensen-2002}.
In this paper we study vertex-critical graphs with respect to the packing chromatic number. The importance of this idea arose during the study of the relationship between the packing chromatic number, chromatic number, and clique number in~\cite{bresar-2018b}. We say that a graph $G$ is {\em packing chromatic vertex-critical}, {\em $\pch$-critical} for short, if $\pch(G-x) < \pch(G)$ holds for every vertex $x$ of $G$. If $\pch(G) = k$, then we also say that $G$ is {\em $k$-$\pch$-critical} (or more simply $k$-critical). To further simplify the terminology we will typically say a {\em critical graph} when considering a packing chromatic  vertex-critical graph.

We proceed as follows. In the next section some definitions are given, known results recalled, and new preliminary results proved. In Section~\ref{sec:vertex-delete} we demonstrate that if $G$ is a $\pch$-critical graph, then the set $\{\pch(G) - \pch(G-x):\ x\in V(G)\}$ can be almost arbitrary. In the subsequent section we characterize $3$-critical graphs and do the same for $4$-critical graphs that contain a cycle of length at least $5$ and not congruent to $0$ modulo $4$. In the case that graphs contain cycles congruent to $0$  modulo $4$, a partial characterization of $4$-critical graphs is given along with some  infinite families of such graphs. In Section~\ref{sec:trees} we concentrate on critical trees and first show that for every integer $k\ge 2$ there exists a $k$-$\pch$-critical tree. In the main result of the section we prove that a $k$-$\pch$-critical caterpillar exists if and only if $k\le 7$. We follow with a section on Cartesian products which turned out to yield many critical graphs. For instance, $K_{1,3}\cp P_3$ is $6$-$\pch$-critical even though neither $K_{1,3}$ nor $P_3$ is critical. In the main result of the section we prove that if $G$ and $H$ are vertex-transitive graphs and $\diam(G) + \diam(H) \le \pch(G)$, then $G\cp H$ is $\pch$-critical. We conclude the paper with several open problems and directions for future investigation.

\section{Preliminaries}

A {\em support vertex} of a graph is a vertex adjacent to a vertex of degree $1$, that is, to a {\em leaf}. By the {\em distance} between two vertices of a graph we mean the usual shortest-path distance. The {\rm diameter}, $\diam(G)$, of a graph $G$ is the maximum distance between all pairs of vertices of $G$. If $H$ is a subgraph of a graph $G$, then for distinct vertices $u$ and $v$ in $H$, the distance between $u$ and $v$ in $H$ is at least
the distance between $u$ and $v$ in $G$.  Hence, if $W$ is an $i$-packing in $G$, it follows that $W \cap V(H)$ is an $i$-packing in $H$.
This gives the following obvious but useful lemma.

\begin{lemma} \label{lem:subgraph}
If $H$ is a subgraph of a graph $G$, then $\pch(G) \ge \pch(H)$.
\end{lemma}

In any packing coloring of a graph with diameter $d$, each color $d$ or larger can be assigned to at most one vertex.  The following is a direct consequence of this fact.

\begin{obser} \label{obs:largecolors}
Suppose $G$ is a graph with $\diam(G)=d \le k=\pch(G)$.  If $c:V(G)\longrightarrow [k]$ is a packing coloring and $x$ is a vertex of $G$ such that $d \le c(x) \le k$, then $\pch(G-x)<\pch(G)$.
\end{obser}

The following result from the seminal paper will be used several times.

\begin{proposition} {\rm (\cite[Proposition 2.1]{goddard-2008})}
\label{prp:indep}
If $G$ is a graph with order $n(G)$, then $\pch(G)\le n(G) - \alpha(G) + 1$, with equality if $\diam(G) = 2$.
\end{proposition}

If $x$ is a vertex of a graph $G$, then the difference between $\pch(G) - \pch(G-x)$ can be made arbitrary large, see~\cite{bkrw-2017a}. The situation is different for leafs as the next observation asserts.

\begin{lemma}
\label{lem:remove-leaf}
If $x$ is a leaf of a graph $G$, then $\pch(G) - 1 \le \pch(G-x) \le \pch(G)$.
\end{lemma}

\proof
Clearly, $\pch(G-x) \le \pch(G)$. Suppose now that $\pch(G-x) = k$. Then using an optimal packing coloring of $G-x$ and using color $k+1$ for the vertex $x$ in $G$ we infer that $\pch(G)\le k+1$. Hence $\pch(G) \le k+1 = \pch(G-x) + 1$.
\qed

We next list some general properties of critical graphs.

\begin{lemma}
\label{lem:connected}
If $G$ is a $\pch$-critical graph, then $G$ is connected.
\end{lemma}

\proof
Suppose $G$ has components $A_1,\ldots,A_r$.  Since $\pch(G)=\max\{\pch(A_i)\}$, there exists a component, say $A_j$, that has packing chromatic number $\pch(G)$. If $r \ge 2$ and $x$ is a vertex in $G-A_j$, then it is clear that $\pch(G-x)=\pch(A_j)=\pch(G)$, which contradicts the assumption that $G$ is $\pch$-critical. Therefore, $r=1$ and hence $G$ is connected.
\qed

If a graph is $k$-critical and some of its edges can be removed without decreasing the packing chromatic number,
then this edge-reduced graph is also $k$-critical as the next result shows.

\begin{proposition} \label{prop:deleteedges}
Let $k$ be a positive integer and let $G$ be a $k$-critical graph.  If $F$ is any subset of the edges of $G$
such that $\pch(G-F)=k$, then $G-F$ is also $k$-critical.
\end{proposition}
\proof
Let $x$ be any vertex in a $k$-critical graph $G$ and let $F$ be any set of edges of $G$ such that $\pch(G-F)=k$.
If $F'$ is that subset of edges in $F$ that are not incident with $x$, then $(G-F)-x=(G-x)-F'$, and $(G-x)-F'$ is
a subgraph of $G-x$.  Since $G$ is $k$-critical, it now follows from Lemma~\ref{lem:subgraph} that
\[\pch((G-F)-x)=\pch((G-x)-F')\le \pch(G-x)<k\,,\]
and therefore $G-F$ is $k$-critical.
\qed

The following corollary is a special case of Proposition~\ref{prop:deleteedges}, but we state it anyway for its later use.

\begin{corollary} \label{cor:addoneedge}
Let $H$ be a graph with $\pch(H)=k$ and let $u$ and $v$ be nonadjacent vertices in $H$.  If $H$ is not $k$-critical, then
$H+uv$ is not $k$-critical.
\end{corollary}

Finally, knowing the packing chromatic number of paths and cycles, we get that $C_n$ is critical if and only if $n\in \{3, 4\}$, or $n\ge 5$ and $n\not\equiv 0 \pmod 4$. It is also easy to verify that if $n\ge 2$, then $K_{n,n}$ is an $(n+1)$-$\pch$-critical graph.

\section{Vertex-deleted subgraphs of $\pch$-critical graphs}
\label{sec:vertex-delete}

In this section we show that the condition $\pch(G-x) < \pch(G)$ for $G$ to be $\pch$-critical cannot be replaced with  $\pch(G-x) = \pch(G) -1$ (as it is the case with the chromatic critical graphs). Moreover, we show that if $G$ is a $\pch$-critical graph, then the set of differences
$$\Delta_{\pch}(G) = \{\pch(G) - \pch(G-x):\ x\in V(G)\}$$
can be almost arbitrary. Before stating the main result, we consider three particular examples.

Clearly, if $G$ is vertex-transitive and critical, then $|\Delta_{\pch}(G)| = 1$. For instance, $\pch(K_{r,r}-w)=\pch(K_{r,r})-1$ for every $w$ and every $r\ge 1$, hence $\Delta_{\pch}(K_{r,r}) = \{1\}$.

Consider next the Petersen graph $P$. Since $\diam(P) = 2$, Proposition~\ref{prp:indep} implies that $\pch(P) = 7$ and in $P-x$ there are three vertices pairwise at distance $3$ that can be colored with $2$. From here a $5$-$\pch$-coloring of $P-x$ is easy to obtain. In conclusion,  $\pch(P-x)=\pch(P)-2$ for every vertex $x$ and thus $\Delta_{\pch}(P) = \{2\}$.

Next, let $2 \le r\le s$, and let $G$ be the graph constructed from identifying a vertex from $K_r$ with a vertex from $K_s$.  The packing chromatic number of $G$ is $r+s-2$.  If $w$ is the vertex shared by the two complete graphs, then $\pch(G-w)=s-1$.  For any other vertex $x$ in $G$, we have $\pch(G-x)=r+s-3$. Hence $\Delta_{\pch}(G) = \{1, r-1\}$.

\begin{theorem}
\label{thm:many-values}
Let $S =\{1,s_1,\ldots,s_r\}$ be a set of positive integers, where $r \ge 1$. If  for every $i\in [r]$ we have $\displaystyle{\sum_{j=1, j\ne i}^r s_j\ge s_i-1}$, then there exists a $\pch$-critical graph $G$ such that $\Delta_{\pch}(G) = S$.
\end{theorem}

\proof
Suppose first that $r\ge 2$. Let $V(K_r) = \{x_1, \ldots, x_r\}$. Let $G(s_1, \ldots, s_r)$ be the graph obtained from $K_r$ such that for every $i\in [r]$, a vertex of a complete graph $X_i$ of order $s_i + 2$ is identified with $x_i$. See Fig.~\ref{fig:G(1,2,4)} for an example of this construction.

\begin{figure}[ht!]
\tikzstyle{every node}=[circle, draw, fill=black!0, inner sep=0pt,minimum width=.2cm]
\begin{center}
\begin{tikzpicture}[thick,scale=.7]
  \draw(0,0) { 
    +(2.50,2.00) -- +(1.25,1.00)
    +(1.25,1.00) -- +(3.75,1.00)
    +(3.75,1.00) -- +(2.50,2.00)
    +(2.08,3.00) -- +(2.92,3.00)
    +(2.92,3.00) -- +(2.50,2.00)
    +(2.08,3.00) -- +(2.50,2.00)
    +(0.00,1.50) -- +(0.00,0.50)
    +(0.00,0.50) -- +(1.25,1.00)
    +(0.00,1.50) -- +(1.25,1.00)
    +(0.42,1.00) -- +(1.25,1.00)
    +(0.42,1.00) -- +(0.00,1.50)
    +(0.42,1.00) -- +(0.00,0.50)
    +(3.75,1.00) -- +(4.58,2.00)
    +(4.58,2.00) -- +(5.83,2.00)
    +(5.83,2.00) -- +(6.67,1.00)
    +(6.67,1.00) -- +(5.83,0.00)
    +(3.75,1.00) -- +(4.58,0.00)
    +(4.58,0.00) -- +(5.83,0.00)
    +(4.58,2.00) -- +(6.67,1.00)
    +(4.58,2.00) -- +(5.83,0.00)
    +(4.58,2.00) -- +(4.58,0.00)
    +(5.83,2.00) -- +(3.75,1.00)
    +(5.83,2.00) -- +(4.58,0.00)
    +(5.83,2.00) -- +(5.83,0.00)
    +(6.67,1.00) -- +(4.58,0.00)
    +(6.67,1.00) -- +(3.75,1.00)
    +(5.83,0.00) -- +(3.75,1.00)
    +(1.25,1.00) node{}
    +(2.50,2.00) node{}
    +(3.75,1.00) node{}
    +(2.08,3.00) node{}
    +(2.92,3.00) node{}
    +(0.00,0.50) node{}
    +(0.00,1.50) node{}
    +(0.42,1.00) node{}
    +(4.58,2.00) node{}
    +(5.83,2.00) node{}
    +(6.67,1.00) node{}
    +(5.83,0.00) node{}
    +(4.58,0.00) node{}
    +(7.35,1) node[rectangle, draw=white!0, fill=white!100]{$X_3$}
    +(2.5,3.5) node[rectangle, draw=white!0, fill=white!100]{$X_1$}
    +(-.65,1) node[rectangle, draw=white!0, fill=white!100]{$X_2$}
  };
\end{tikzpicture}
\end{center}
\caption{The graph $G(1,2,4)$.} \label{fig:G(1,2,4)}
\end{figure}
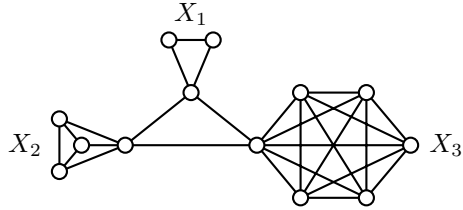

To simplify the notation, set $G = G(s_1, \ldots, s_r)$ in the sequel. Note first that $$|V(G)| = \sum_{i=1}^r |V(X_i)| = \sum_{i=1}^r (s_i + 2) =  \sum_{i=1}^r s_i + 2r\,.$$
If $c$ is a $\pch$-coloring of $G$, then the vertices of $X_i$, $i\in [r]$, receive different colors. In particular, $|c^{-1}(1)| \le r$ and $|c^{-1}(2)| \le r$. Moreover, since $\diam(G) = 3$, $|c^{-1}(\ell)| \le 1$ for any $\ell \ge 3$. Since $s_i\ge 2$ and so $s_i+2 \ge 4$, $i\in [r]$, in each $X_i$ colors $1$ and $2$ can be used. It follows  that
\begin{equation}
\label{eq:pch-number-of-G}
\pch(G) = 2 + (|V(G)| - 2r) = 2 + \sum_{i=1}^r s_i\,.
\end{equation}
Since $r\ge 2$, for at least one $s_i$ we have $s_i\ge 3$. Assume without loss of generality that $s_1\ge 3$. Let $u\in X_1$ be an arbitrary vertex different from $x_1$. Then $G-u$ is isomorphic to $G(s_1-1, s_2, \ldots, s_r)$ (where it is possible that $s_1-1 = s_i$ for some $i\ge 2$). Then by~\eqref{eq:pch-number-of-G} we get $\pch(G-u) = 2 + (s_1-1) + \sum_{i=2}^r s_i =  \sum_{i=1}^r s_i + 1 = \pch(G) - 1$. This shows that $1\in \Delta_{\pch}(G)$.

Consider next the vertex deleted subgraph $G-x_i$, $i\in [r]$. Since $x_i$ is a cut vertex, we infer, having in mind~\eqref{eq:pch-number-of-G}, that
\begin{eqnarray*}
\pch(G-x_i) & = & \max \{ \pch(K_{s_i+1}), \pch(G(s_1,\ldots, s_{i-1}, s_{i+1}, \ldots, s_r)) \} \\
& = & \max \{ s_i+1, 2 + \sum_{j=1, j\ne i}^r s_j \} \\
& = & 2 + \sum_{j=1, j\ne i}^r s_j\,,
\end{eqnarray*}
where the last equality follows since we have assumed that $\displaystyle{\sum_{j=1, j\ne i}^r s_j\ge s_i-1}$. From here and~\eqref{eq:pch-number-of-G} we get that $\pch(G) - \pch(G-x_i) = s_i$, that is, $s_i\in \Delta_{\pch}(G)$ for $i\in [r]$. This proves the theorem for $r\ge 2$.

Let now $S = \{1, s\}$, where $s\ge 2$. Let $G$ be the graph obtained from two disjoint copies of $K_{s+1}$ by identifying a vertex from one copy with a vertex from the other copy. Let $x$ be the identified vertex. Then $|V(G)| = 2s+1$ and $\pch(G) = 1 + (2s+1 - 2) = 2s$. If $u$ is a vertex of $G$ different from $x$, then $\pch(G-u) = 1 + (2s-2) = 2s -1$, so that $\pch(G) - \pch(G-u) = 1$. Finally $\pch(G-x) = \pch(K_s) = s$, so that $\pch(G) - \pch(G-x) = s$ and we are done.
\qed

\section{On $3$-$\pch$-critical and $4$-$\pch$-critical graphs}
\label{sec:3-and-4-critical}

Clearly, $K_2$ is the unique $2$-$\pch$-critical graph. To characterize the $3$-$\pch$-critical graphs we recall from~\cite[Proposition 3.1]{goddard-2008} that if $G$ is a connected graph, then $\pch(G) = 2$ if and only if $G$ is a nontrivial star.

\begin{proposition}
\label{prp:3-critical}
A graph $G$ is $3$-$\pch$-critical if and only if $G\in \{C_3, P_4, C_4\}$.
\end{proposition}

\proof
Let $G$ be a $3$-$\pch$-critical graph.  Then $G$ is connected by Lemma~\ref{lem:connected}.  Clearly, $|V(G)|\ge 3$. If $|V(G)| = 3$, then
$G\in \{C_3, P_3\}$, giving us the $3$-$\pch$-critical graph $C_3$. Assume in the rest that $|V(G)|\ge 4$ and let $u$ be an arbitrary vertex of $G$. Then $\pch(G-u) = 2$ and hence by the above remark, $G-u$ is a disjoint union of stars. If at least two of these stars contain an edge, then $G$ contains a $P_5$ and hence cannot be $3$-$\pch$-critical. On the other hand, at least one of the stars, say $S_1$, must contain an edge, for otherwise $G$, being connected, would itself be a star. Then $G-u$ contains at most one isolated vertex, for otherwise removing one such vertex from $G$ would yield a graph with $\pch = 3$. To summarize thus far, $G-u$ contains the star $S_1$ with at least one edge and contains at most one isolated vertex.

Suppose first that $G-u = S_1$. As $|V(G)|\ge 4$, $S_1$ has at least two leaves. If $u$ is adjacent to the center of $S_1$, then since $G$ is a not a star, $u$ must be adjacent to at least one leaf of $S_1$. Removing in $G$ another leaf of $S_1$ leaves a graph with $\pch = 3$ (it contains a $C_3$). So $u$ is not adjacent to the center of $S_1$ and is hence adjacent to at least one leaf in $S_1$. If $S_1$ has at least three leaves, then $G$ again cannot be $3$-$\pch$-critical. In conclusion, $S_1$ has exactly two leaves. If $u$ is adjacent to exactly one of them, we get $P_4$, and if it is adjacent to both, we get $C_4$.

Suppose next that $G-u$ is the disjoint union of $S_1$ and an isolated vertex, say  $w$. If $S_1$ has at least two leaves, then similarly as in the first case we infer that $G$ is not $3$-$\pch$-critical.  And if $S_1 = K_2$, then we either find the $3$-$\pch$-critical $P_4$ or the triangle with a pendant edge which is not $3$-$\pch$-critical.

Since it is clear that each of $C_3$, $P_4$, and  $C_4$ is $3$-$\pch$-critical, the proof is complete.
\qed

We next turn our attention to $4$-critical graphs.

\begin{lemma} \label{lem:propercycle}
If a graph $G$ contains a cycle $C$ whose order is larger than $3$ and is not divisible by $4$ such that $V(G)-V(C) \not=\emptyset$, then $G$ is not $4$-critical.
\end{lemma}

\proof
Suppose that $\pch(G) = 4$, for otherwise there is nothing to be proved. Let $C$ be a cycle in $G$ of order $n>3$ such that $n \not\equiv 0 \pmod{4}$ and $V(G) \not=V(C)$.  For any $x \in V(G-C)$, the cycle $C$ is a subgraph of $G-x$ and thus by Lemma~\ref{lem:subgraph}, $4=\pch(G) \ge \pch(G-x) \ge \pch(C)=4$.  Therefore, $G$ is not $4$-critical.
\qed

For the sake of simplicity we call a cycle in a graph $G$ that does not contain all the vertices of $G$ a {\em proper cycle}.  We are now able to characterize the $4$-$\pch$-critical graphs that contain a cycle of order at least $5$ that is not divisible by $4$.

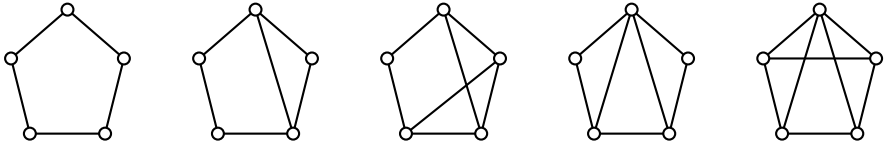
\begin{figure}[ht!]
\begin{center}
\begin{tikzpicture}[scale=.5,style=thick]
\def\vr{4.5pt} 
\path (0,0) coordinate (a1); \path (-0.5,2) coordinate (a2); \path (1,3.3) coordinate (a3); \path (2.5,2) coordinate (a4); \path (2,0) coordinate (a5);
\path (5,0) coordinate (b1); \path (4.5,2) coordinate (b2); \path (6,3.3) coordinate (b3); \path (7.5,2) coordinate (b4); \path (7,0) coordinate (b5);
\path (10,0) coordinate (c1); \path (9.5,2) coordinate (c2); \path (11,3.3) coordinate (c3); \path (12.5,2) coordinate (c4); \path (12,0) coordinate (c5);
\path (15,0) coordinate (d1); \path (14.5,2) coordinate (d2); \path (16,3.3) coordinate (d3); \path (17.5,2) coordinate (d4); \path (17,0) coordinate (d5);
\path (20,0) coordinate (e1); \path (19.5,2) coordinate (e2); \path (21,3.3) coordinate (e3); \path (22.5,2) coordinate (e4); \path (22,0) coordinate (e5);
\draw (a1) -- (a2); \draw (a2) -- (a3);  \draw (a3) -- (a4); \draw (a4) -- (a5); \draw (a5) -- (a1);
\draw (b1) -- (b2); \draw (b2) -- (b3);  \draw (b3) -- (b4); \draw (b4) -- (b5); \draw (b5) -- (b1); \draw (b3) -- (b5);
\draw (c1) -- (c2); \draw (c2) -- (c3);  \draw (c3) -- (c4); \draw (c4) -- (c5); \draw (c5) -- (c1); \draw (c3) -- (c5); \draw (c4) -- (c1);
\draw (d1) -- (d2); \draw (d2) -- (d3);  \draw (d3) -- (d4); \draw (d4) -- (d5); \draw (d5) -- (d1); \draw (d3) -- (d5); \draw (d3) -- (d1);
\draw (e1) -- (e2); \draw (e2) -- (e3);  \draw (e3) -- (e4); \draw (e4) -- (e5); \draw (e5) -- (e1); \draw (e3) -- (e5); \draw (e3) -- (e1); \draw (e2) -- (e4);
\foreach \i in {1,...,5}
{  \draw (a\i)  [fill=white] circle (\vr); \draw (b\i)  [fill=white] circle (\vr); \draw (c\i)  [fill=white] circle (\vr);
\draw (d\i)  [fill=white] circle (\vr); \draw (e\i)  [fill=white] circle (\vr);
 }
\end{tikzpicture}
\end{center}
\caption{$4$-critical graphs obtained from $C_5$ by adding chords.}
\label{fig:order5}
\end{figure}

\begin{theorem} \label{thm:largecycle}
If $G$ is a graph that contains a cycle of length $n\ge 5$, where $n \not\equiv 0 \pmod{4}$, then $G$ is  $4$-critical if and only if one of the following holds.
\begin{itemize}
\item $n=5$ and $G$ is one of the graphs in Fig.~\ref{fig:order5},
\item $n=6$ and $G$ is one of four bipartite graphs obtained by adding some subset of chords to the $6$-cycle,
\item $n\ge 7$ and $G$ is isomorphic to $C_n$.
\end{itemize}
\end{theorem}

\proof
Let $n$ be an integer, $n \ge 5$, such that $n \not\equiv 0 \pmod{4}$.  Note first that $\pch(C_n)=4$ and each vertex-deleted subgraph of $C_n$ is the path $P_{n-1}$,
which has packing chromatic number $3$.  Thus $C_n$ is $4$-$\pch$-critical.  Let $G$ be a $4$-$\pch$-critical graph that contains a cycle $C$ of order $n$.  By
Lemma~\ref{lem:propercycle} it follows that $C$ is not a proper cycle; that is, $C$ is a Hamiltonian cycle.  We let $C$ have vertices $v_0,v_1, \ldots,v_{n-1}$ and edges $v_iv_{i+1}$ for $0 \le i \le n-1$, where the subscripts are computed modulo $n$.  A chord $v_iv_j$ of $C$ is called an \emph{$r$-chord} if the shortest distance between $v_i$ and $v_j$ on $C$ is $r$.  Note that $C$ has $r$-chords for $2 \le r \le \lfloor \frac{n}{2}\rfloor$.

First we let $n=5$.  All the chords of the $5$-cycle are $2$-chords.  One can easily check that the five graphs of Fig.~\ref{fig:order5} are $4$-critical and that
if chords are added in any other way, the resulting graph is not $4$-critical.  Second, let $n=6$.  By Proposition~\ref{prp:indep}, $\pch(K_{3,3})=4$ and since each vertex-deleted subgraph of $K_{3,3}$ is isomorphic to $K_{2,3}$, Proposition~\ref{prp:indep} implies that $\pch(K_{2,3})=3$.  Therefore, $K_{3,3}$ is $4$-$\pch$-critical.  Since $C_6$ is also 4-critical, it follows from Proposition~\ref{prop:deleteedges} that the two subgraphs of $K_{3,3}$ obtained by deleting one or two  independent edges from $K_{3,3}$ are all $4$-critical.  On the other hand, if any $2$-chord is added to a $6$-cycle, the resulting graph contains a proper cycle of order $5$ and is therefore not $4$-critical.

Now, let $n \ge 7$.  Let $G_r$ be the graph formed by adding the $r$-chord $v_0v_r$ to $C$.
We claim that for each $2 \le r \le \lfloor \frac{n}{2}\rfloor$, the graph $G_r$ is not $4$-critical.  It will then follow by  Corollary~\ref{cor:addoneedge} that no supergraph of $G_r$ is $4$-critical.  This will show that $G$ is isomorphic to the $n$-cycle.

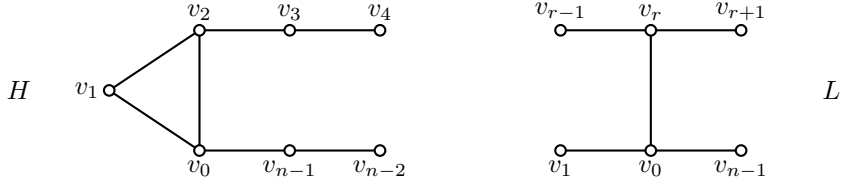
\begin{figure}[ht!]
\begin{center}
\begin{tikzpicture}[scale=0.4,style=thick]
\def\vr{5pt} 
\path (0,0) coordinate (v1); \path (3,-2) coordinate (v0); \path (3,2) coordinate (v2);
\path (6,-2) coordinate (vn-1); \path (6,2) coordinate (v3); \path (9,-2) coordinate (vn-2);
\path (9,2) coordinate (v4);
\path (15,-2) coordinate (w1); \path (15,2) coordinate (wr-1); \path (18,2) coordinate (wr);
\path (18,-2) coordinate (w0); \path (21,2) coordinate (wr+1); \path (21,-2) coordinate (wn-1);
\draw (v0) -- (v1); \draw (v0) -- (v2); \draw (v1) -- (v2); \draw (v2) -- (v3);
\draw (v3) -- (v4); \draw (v0) -- (vn-1); \draw (vn-1) -- (vn-2);
\draw (w0) -- (w1); \draw (w0) -- (wr);  \draw (wr-1) -- (wr);
\draw (wr) -- (wr+1); \draw (w0) -- (wn-1);
%
\foreach \i in {0,1,2,3,4}
{  \draw (v\i)  [fill=white] circle (\vr);
 }
\draw (vn-1)  [fill=white] circle (\vr); \draw (vn-2)  [fill=white] circle (\vr);

\foreach \i in {0,1}
{  \draw (w\i)  [fill=white] circle (\vr);
 }
\draw (wn-1)  [fill=white] circle (\vr); \draw (wr+1)  [fill=white] circle (\vr); \draw (wr)  [fill=white] circle (\vr); \draw (wr-1)  [fill=white] circle (\vr);

\draw[anchor=east] (v1) node{$v_1$}; \draw[anchor=south] (v2) node{$v_2$}; \draw[anchor=south] (v3) node{$v_3$};
\draw[anchor=south] (v4) node{$v_4$}; \draw[anchor=north] (v0) node{$v_0$};
\draw[anchor=north] (vn-1) node{$v_{n-1}$};  \draw[anchor=north] (vn-2) node{$v_{n-2}$};

\draw[anchor=north] (w1) node{$v_1$}; \draw[anchor=south] (wr-1) node{$v_{r-1}$}; \draw[anchor=south] (wr) node{$v_r$};
\draw[anchor=south] (wr+1) node{$v_{r+1}$}; \draw[anchor=north] (w0) node{$v_0$};
\draw[anchor=north] (wn-1) node{$v_{n-1}$};

\draw (-3,0) node{$H$};
\draw (24,0) node{$L$};

\end{tikzpicture}
\end{center}
\caption{Subgraphs $H$ and $L$.} \label{fig:2-chord}
\end{figure}

We consider first the case $r=2$.  The graph $G_2$ contains proper cycles of order $3$ and $n-1$.  If $n-1 \not\equiv 0 \pmod 4$, then  by Lemma~\ref{lem:propercycle} we conclude  that $G_2$ is not $4$-critical.  Therefore, we assume that $n\equiv 1 \pmod{4}$.
In this case $n \ge 9$ and we infer that $G_r$ contains a vertex $x$ such that $G_r-x$ has a subgraph isomorphic to the graph $H$ in Fig.~\ref{fig:2-chord}.
We claim that $\pch(H)\ge 4$.  Let $c:V(H) \to [k]$ be a $k$-packing coloring.  If
$c(v_0)=1$, then $c(v_1)$, $c(v_2)$ and $c(v_{n-1})$ are three distinct colors, which implies $k \ge 4$.  Similarly, $c(v_2)=1$ implies $k \ge 4$.  Hence, $k\ge 4$ unless $\{c(v_0),c(v_2)\}=\{2,3\}$.  However, if
$\{c(v_0),c(v_2)\}=\{2,3\}$, then it is easy to see that one of $c(v_1)$, $c(v_3)$, $c(v_4)$, $c(v_{n-1})$ or
$c(v_{n-2})$ is at least 4.  Therefore, $\pch(H) \ge 4$ and by Lemma~\ref{lem:subgraph} it follows that $G_2$ is not $4$-critical.

Now assume that $r \ge 3$.  Recall that $n \ge 7$.  The graph $G_r$ has proper cycles of orders $r+1$ and $n-r+1$.  Note that
$r+1 > 3$ and $n-r+1 >3$.  A straightforward computation shows that either $r+1 \not\equiv 0 \pmod{4}$ or
$n-r+1 \not\equiv 0 \pmod{4}$ except when $n\equiv 2 \pmod{4}$ and $r \equiv 3 \pmod{4}$.  Therefore, if $n\not \equiv 2 \pmod{4}$ or $r \not \equiv 3 \pmod{4}$,
then it follows by  Lemma~\ref{lem:propercycle} that $G_r$ is not $4$-critical.
Thus we assume that $n\equiv 2 \pmod{4}$ and $r \equiv 3 \pmod{4}$. Suppose that $G_r$ is $4$-critical and let $c: V(G_r-v_{r+2}) \to [3]$  be a $3$-packing
coloring of $G_r-v_{r+2}$.   Let $L$ be the subgraph of $G_r-v_{r+2}$ shown in Fig.~\ref{fig:2-chord}. If $c(v_0)=1$, then $c(v_1)$, $c(v_{n-1})$ and $c(v_r)$
must all be different, which is a contradiction.   Therefore, $c(v_0)\neq 1$, and similarly $c(v_r)\neq 1$.  Without loss of generality we assume
that $c(v_0)=3$ and $c(v_r)=2$.  This implies that $c(v_1)=1$ and $c(v_{r-1})=1$.  If $r=3$, we arrive at a contradiction since then $v_1$ and $v_{r-1}$
are adjacent.  For $r>3$ the vertex $v_{r-2}$ is distance $2$ from $v_r$ and distance $3$ from $v_0$, which is again a contradiction.  This implies
that $\pch(G_r-v_{r+2}) \ge 4$ and hence $G_r$ is not 4-critical.

Therefore, for $n\ge 7$ such that  $n \not \equiv 0 \pmod{4}$ we have shown that if a $4$-$\pch$-critical graph $G$ contains a cycle of order $n$, then
$G=C_n$. This completes the proof of the theorem.
\qed

There remains the case when the longest cycle of a 4-packing chromatic graph has order $n$ congruent to 0 modulo $4$.  If such a
graph is $4$-critical and $n \ge 8$, then none of its cycles of order eight or more can have chords since any such chord creates a proper cycle whose order is at least $5$ and is not
congruent to 0 modulo $4$. For a partial characterization of $4$-critical graphs all of whose cycles have order congruent to 0 modulo $4$ we first prove the following.

\begin{lemma}
\label{lem:even-path-with-2-leaves}
Let $G_{2k}$, $k\ge 3$, be the graph obtained from the path $P_{2k}$ by attaching a leaf to each of the vertices at distance $2$ from an endvertex. Then $\pch(G_{2k}) = 4$.
\end{lemma}

\proof
Let the vertices of $P_{2k}$ be $v_1,\ldots, v_{2k}$, and let $v_3'$ and $v_{2k-2}'$ be the leaves attached to $v_3$ and $v_{2k-2}$, respectively. Coloring the vertices of $P_{2k}$ with the pattern $1,2,1,3,1,2,1,3,\ldots$, and the vertices $v_3'$ and $v_{2k-2}'$ with $4$ and $1$, respectively, we get $\pch(G_{2k}) \le 4$.
For the rest of the proof we first state the following observation.

\medskip\noindent
{\bf Fact}. In any $3$-$\pch$-coloring of a path, if $x$ and $y$ are adjacent vertices both of which have a distance at least $2$ from the ends of the path, then at least one of $x$ and $y$ must be colored $1$.

Suppose now that $G_{2k}$ admits a $3$-packing coloring $c$. Then $c(v_3)\ne 1$ and $c(v_{2k-2})\ne 1$. Suppose that $c(v_3) = 2$. By the above Fact $c(v_4) \ne 3$, that is, $c(v_4) = 1$. If $k=3$, this is a contradiction, so let $k > 3$. Then $c(v_5) = 3$. Repeatedly using the Fact we deduce that $c(v_{2k-2}) = 1$ must hold, a contradiction again. Therefore $c(v_3) = 3$. Then we first infer that  either $c(v_2) = 2$ and $c(v_1) = 1$, or $c(v_2) = 1$ and $c(v_1) = 2$. In the first case $c(v_4) = 1$ clearly holds. In the second case, assuming that $c(v_4) = 2$, we get that $v_6$ cannot be properly colored. Hence in either case, $c(v_4) = 1$. But then, using the Fact as many times as required it follows again that $c(v_{2k-2}) = 1$. Hence $c$ does not exist and so $\pch(G_{2k}) \ge 4$.
\qed

Let ${\cal C}$ be the class of graphs that contain exactly one cycle and have an arbitrary number of attached leaves to each of the vertices of the cycle. Recall that the {\em net graph} is the graph formed by attaching a single leaf to each vertex of a $K_3$. Then the announced partial characterization of  $4$-critical graphs reads as follows.

\begin{theorem}
\label{thm:C-cycles}
A graph $G\in {\cal C}$ is $4$-$\pch$-critical if and only if $G$ is one of the following graphs:
\begin{itemize}
\item $G=C_n$,  $n\ge 5$, $n\not\equiv 0 \pmod 4$,
\item $G$ is the net graph,
\item $G$ is obtained by attaching a single leaf to two adjacent vertices of $C_4$.
\item $G$ is obtained by attaching a single leaf to two  vertices at distance $3$ on $C_8$.
\end{itemize}
\end{theorem}

\proof
Let $G\in {\cal C}$, let $C$ be the unique cycle in $G$, and assume $C$ has order $n$. If $n\not\equiv 0 \pmod 4$ and $n\ge 5$, then $G$ is not $4$-critical by Lemma~\ref{lem:propercycle} if $G$ has at least one leaf. On the other hand it is straightforward to see that $C_n$ ($n\not\equiv 0 \pmod 4$ and $n\ge 5$) is $4$-critical. If $n=3$, then it is also equally straightforward to verify that the net graph is the only graph from ${\cal C}$ that contains a $3$-cycle and is $4$-critical. Hence we may assume in the rest that $n = 4k$ for some $k\ge 1$.

Suppose that $G$ is $4$-$\pch$-critical. We first claim $G$ has at most one leaf attached to each of the vertices of $C$. Suppose on the contrary that a vertex $u$ of $C$ is adjacent to two leaves $x$ and $y$. Since $G$ is $4$-$\pch$-critical, $\pch(G-x) = 3$. In a $3$-coloring of $G-x$, the vertex $u$ is not colored $1$. But then the $3$-coloring of $G-x$ can be extended to a $3$-coloring of $G$ by coloring $x$ with $1$, a contradiction which proves the claim.

Suppose first that each pair of support vertices of $G$ is at even distance. Then we can color $C$ with the pattern $1,2,1,3,1,2,1,3,\ldots$ such that each support vertex receives color $2$ or $3$. But then $\pch(G) = 3$. Similarly, if $G$ has only one support vertex, then we also have $\pch(G) = 3$. This implies that $G$ is not $4$-critical, which is a contradiction.

Therefore $G$ has at least two support vertices at odd distance. Let $u$ and $v$ be two support vertices that are closest possible on $C$ and let $d_G(u,v) = \ell$. Note that $\ell \le 2k-1$. Suppose that $k\ge 3$ and let $x$ be a vertex that is not on a $u,v$-shortest path and $d_G(x,u)\ge 3$ and $d_G(x,v)\ge 3$. As $k\ge 3$, such a vertex $x$ exists. Then $G-x$ contains a subgraph isomorphic to the graph $G_{\ell+5}$ from Lemma~\ref{lem:even-path-with-2-leaves} which in turn implies that $G$ is not $4$-critical. Assume next that $k=2$, that is, $C=C_8$. If $\ell = 1$, the same argument implies that $G$ is not $4$-critical. If $\ell = 3$, then Lemma~\ref{lem:even-path-with-2-leaves} yields $\pch(G) = 4$. On the other hand it is straightforward to verify that $\pch(G-x) = 3$ for every $x\in V(G)$. Hence $G$ is $4$-critical. Finally, if $k=1$, then we only need to check two graphs. Among them the one that is obtained by attaching a single leaf to two adjacent vertices of $C_4$ is $4$-critical.
\qed

To conclude the section we present another infinite family of $4$-$\pch$-critical graphs. Let $H_{2k+1}$, $k\ge 0$, be the graph obtained from the disjoint union of two $4$-cycles by connecting a vertex of one $4$-cycle with a vertex of the other $4$-cycle with a path of length $2k+1$. Then applying Lemma~\ref{lem:even-path-with-2-leaves} we see that $\pch(H_{2k+1}) = 4$. If $x$ is a diametrical vertex of $H_{2k+1}$, then distinguishing the cases when $2k+1 \equiv 3 \pmod 4$ and $2k+1 \equiv 1 \pmod 4$ we see that $\pch(H_{2k+1}-x) = 3$, see Fig.~\ref{f:G}.

\begin{figure}[ht!]
\begin{center}
\begin{tikzpicture}[scale=1.0,style=thick]
\def\vr{2.5pt} 
\path (2.00,3.33) coordinate (v1);
\path (3.00,3.33) coordinate (v2);
\path (4.00,3.33) coordinate (v3);
\path (5.00,3.33) coordinate (v4);
\path (5.50,4.00) coordinate (v5);
\path (6.00,3.33) coordinate (v6);
\path (5.50,2.67) coordinate (v7);
\path (1.50,4.00) coordinate (v8);
\path (1.00,3.33) coordinate (v99);
\path (1.50,2.67) coordinate (v10);
\path (2.00,0.67) coordinate (v11);
\path (1.00,0.67) coordinate (v12);
\path (3.00,0.67) coordinate (v13);
\path (4.00,0.67) coordinate (v14);
\path (5.00,0.67) coordinate (v15);
\path (6.00,0.67) coordinate (v16);
\path (6.50,1.33) coordinate (v17);
\path (7.00,0.67) coordinate (v18);
\path (6.50,0.00) coordinate (v19);
\path (0.50,1.33) coordinate (v20);
\path (0.00,0.67) coordinate (v31);
\path (0.50,0.00) coordinate (v22);

\draw (v8) -- (v1); \draw (v1) -- (v10);  \draw (v1) -- (v2); \draw (v2) -- (v3);
\draw (v3) -- (v4); \draw (v4) -- (v5); \draw (v4) -- (v7); \draw (v5) -- (v6); \draw (v6) -- (v7);
\draw (v20) -- (v12); \draw (v12) -- (v22);  \draw (v12) -- (v11); \draw (v11) -- (v13);
\draw (v13) -- (v14); \draw (v14) -- (v15); \draw (v15) -- (v16); \draw (v16) -- (v17); \draw (v17) -- (v18); \draw (v18) -- (v19);
\draw (v16) -- (v19);

\foreach \i in {1,...,8}
{  \draw (v\i)  [fill=white] circle (\vr);
 }
\foreach \i in {10,...,20}
{  \draw (v\i)  [fill=white] circle (\vr);
 }
\draw (v22)  [fill=white] circle (\vr);
\draw[anchor=east] (v8) node{$1$}; \draw[anchor=east] (v22) node{$1$}; \draw[anchor=east] (v20) node{$1$}; \draw[anchor=east] (v10) node{$1$};
\draw[anchor=north] (v1) node{$2$}; \draw[anchor=north] (v2) node{$3$}; \draw[anchor=north] (v3) node{$1$}; \draw[anchor=north] (v4) node{$2$};
\draw[anchor=north] (v7) node{$1$}; \draw[anchor=north] (v12) node{$2$}; \draw[anchor=north] (v11) node{$3$}; \draw[anchor=north] (v13) node{$1$};
\draw[anchor=south] (v5) node{$1$}; \draw[anchor=south] (v17) node{$1$}; \draw[anchor=west] (v6) node{$3$}; \draw[anchor=west] (v18) node{$2$};
\draw[anchor=north] (v14) node{$2$}; \draw[anchor=north] (v15) node{$1$}; \draw[anchor=north] (v16) node{$3$}; \draw[anchor=north] (v19) node{$1$};
\end{tikzpicture}
\end{center}
\caption{$3$-$\pch$-colorings of $H_3-x$ and of $H_5-x$.} \label{f:G}
\end{figure}

For all the other vertices $y$ of $H_{2k+1}$ it is easy to see that $\pch(H_{2k+1}-y) = 3$. Hence $H_{2k+1}$ is
$4$-critical for every $k\ge 0$. On the other hand, if  in the construction of $H_{2k+1}$ one or both $4$-cycles are replaced with some longer cycle of order congruent to $0$ modulo $4$, then the obtained graph is not $4$-critical as can be again deduced from Lemma~\ref{lem:even-path-with-2-leaves}.

\section{On critical trees}
\label{sec:trees}

We begin this section with the following interesting fact.

\begin{proposition}
\label{prop:trees-for-all-k}
If $k\ge 2$, then there exists a $k$-$\pch$-critical tree.
\end{proposition}

\proof
\cite{sloper-2004} proved that the infinite $4$-regular tree has no finite packing coloring. Hence there exists a finite tree $T$ with $\pch(T) = k$. Indeed, by Sloper's result there exists a finite tree with packing chromatic number $s\ge k$. If $s>k$, then we repeatedly use Lemma~\ref{lem:remove-leaf} to arrive at the desired tree. Now, if $T$ is critical, we are done. Otherwise there exists a vertex $x$ of $T$ such that $\pch(T-x) = \pch(T)$. It follows that $T-x$ contains a component $T'$ with $\pch(T') = k$. If $T'$ is critical, we are done, otherwise just continue the process until a $k$-$\pch$-critical tree is found.
\qed

\begin{figure}[th!]
\tikzstyle{every node}=[circle, draw, fill=black!0, inner sep=0pt,minimum width=.2cm]
\begin{center}
\begin{tikzpicture}[thick,scale=1.5]
  \draw(0,0) { 
    +(0.00,3.33) -- +(0.50,3.33)
    +(0.50,3.33) -- +(0.50,2.67)
    +(0.00,3.33) -- +(0.00,2.67)
    +(1.50,2.67) -- +(1.50,3.33)
    +(1.50,3.33) -- +(2.00,3.33)
    +(2.00,3.33) -- +(2.00,2.67)
    +(2.00,3.33) -- +(2.50,3.33)
    +(2.50,3.33) -- +(2.50,2.67)
    +(2.50,3.33) -- +(3.00,3.33)
    +(3.00,3.33) -- +(3.00,2.67)
    +(4.50,3.33) -- +(5.00,3.33)
    +(5.00,3.33) -- +(5.50,3.33)
    +(5.50,3.33) -- +(5.50,2.67)
    +(5.00,3.33) -- +(5.00,2.67)
    +(4.50,3.33) -- +(4.50,2.67)
    +(0.00,0.67) -- +(0.50,0.67)
    +(0.50,0.67) -- +(1.00,0.67)
    +(1.00,0.67) -- +(1.50,0.67)
    +(1.50,0.67) -- +(2.00,0.67)
    +(2.00,0.67) -- +(2.50,0.67)
    +(2.50,0.67) -- +(2.50,0.00)
    +(2.00,0.67) -- +(2.00,0.00)
    +(0.50,0.67) -- +(0.50,0.00)
    +(0.00,0.67) -- +(0.00,0.00)
    +(4.00,1.33) -- +(4.00,0.67)
    +(4.00,0.67) -- +(4.00,0.00)
    +(4.00,0.67) -- +(4.50,0.67)
    +(4.50,0.67) -- +(5.00,0.67)
    +(5.00,0.67) -- +(5.50,0.67)
    +(5.50,0.67) -- +(6.00,0.67)
    +(6.00,1.33) -- +(6.00,0.67)
    +(6.00,0.67) -- +(6.00,0.00)
    +(5.50,1.33) -- +(5.50,0.67)
    +(5.50,0.67) -- +(5.50,0.00)
    +(5.00,1.33) -- +(5.00,0.67)
    +(5.00,0.67) -- +(5.00,0.00)
    +(4.50,1.33) -- +(4.50,0.67)
    +(4.50,0.67) -- +(4.50,0.00)
    +(4.50,4.00) -- +(4.50,3.33)
    +(5.50,4.00) -- +(5.50,3.33)
    +(0.00,3.33) node{}
    +(0.50,3.33) node{}
    +(0.00,2.67) node{}
    +(0.50,2.67) node{}
    +(1.50,3.33) node{}
    +(2.00,3.33) node{}
    +(2.50,3.33) node{}
    +(3.00,3.33) node{}
    +(4.50,3.33) node{}
    +(5.00,3.33) node{}
    +(5.50,3.33) node{}
    +(5.50,2.67) node{}
    +(5.00,2.67) node{}
    +(4.50,2.67) node{}
    +(1.50,2.67) node{}
    +(2.00,2.67) node{}
    +(2.50,2.67) node{}
    +(3.00,2.67) node{}
    +(0.00,0.67) node{}
    +(0.50,0.67) node{}
    +(2.00,0.67) node{}
    +(2.50,0.67) node{}
    +(4.00,1.33) node{}
    +(4.50,1.33) node{}
    +(5.00,1.33) node{}
    +(5.50,1.33) node{}
    +(6.00,1.33) node{}
    +(1.00,0.67) node{}
    +(1.50,0.67) node{}
    +(0.00,0.00) node{}
    +(0.50,0.00) node{}
    +(2.00,0.00) node{}
    +(2.50,0.00) node{}
    +(4.00,0.67) node{}
    +(4.50,0.67) node{}
    +(5.00,0.67) node{}
    +(5.50,0.67) node{}
    +(6.00,0.67) node{}
    +(4.00,0.00) node{}
    +(4.50,0.00) node{}
    +(5.00,0.00) node{}
    +(5.50,0.00) node{}
    +(6.00,0.00) node{}
    +(5.50,4.00) node{}
    +(4.50,4.00) node{}
    +(0.2,2.2) node[rectangle, draw=white!0, fill=white!100]{$\pch=3$}
    +(1.20,-0.5) node[rectangle, draw=white!0, fill=white!100]{$\pch=4$}
    +(2.2,2.2) node[rectangle, draw=white!0, fill=white!100]{$\pch=4$}   
    +(5,2.2) node[rectangle, draw=white!0, fill=white!100]{$\pch=4$}
    +(5.0,-0.5) node[rectangle, draw=white!0, fill=white!100]{$\pch=5$}
  };
\end{tikzpicture}
\end{center}
\caption{Some $k$-$\pch$-critical caterpillars.}
\label{fig:small-critical-trees}
\end{figure}
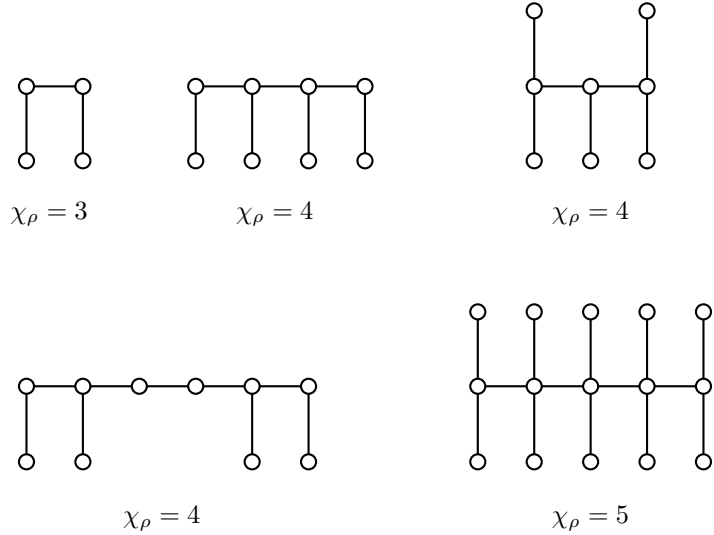

In Fig.~\ref{fig:small-critical-trees} five $k$-$\pch$-critical trees are given, where $k\in \{3,4,5\}$.

Note that all the critical trees from Fig.~\ref{fig:small-critical-trees} are caterpillars. Hence it is natural to wonder for which $k$ there exist $k$-$\pch$-critical caterpillars.

\begin{theorem}
\label{thm:caterpillars}
A $k$-$\pch$-critical caterpillar exists if and only if $k\le 7$.  \end{theorem}

\proof
Let $P_{\infty}$ denote the one-way infinite path with vertices $v_1,v_2,\ldots$ such that $v_1$ has degree 1. By repeating the coloring pattern
$2, 4, 3, 2, 5, 6, 2, 4, 3, 2, 5, 7$ (mentioned in~\cite{sloper-2004}) starting with $v_1$ one can see that $P_{\infty}$ admits a $(2,3,4,5,6,7)$-coloring. If $T$ is a caterpillar, then using the described $(2,3,4,5,6,7)$-coloring for its spine and color $1$ for each of its leaves attached to the spine we infer that $\pch(T)\le 7$ holds for an arbitrary caterpillar $T$.

\cite{sloper-2004} also mentioned that no $(2,3,4,5,6)$-coloring of $P_n$ exists if $n\ge 35$. (We have independently verified this fact using a backtracking algorithm.) Hence, if $n$ is such, then attaching to every vertex of $P_n$ precisely (or more) $6$ leaves yields a caterpillar $T$ with $\pch(T) = 7$. Iteratively applying Lemma~\ref{lem:remove-leaf} to the leaves of $T$ we find a $7$-$\pch$-critical caterpillar. Continuing in this manner we then see that there also exists a $k$-$\pch$-critical caterpillar for every $k\le 6$.
\qed

Note that the proof of Theorem~\ref{thm:caterpillars} is not constructive. On the other hand, Fig.~\ref{fig:small-critical-trees} shows $k$-$\pch$-critical caterpillars for $k\le 5$. In the next result we construct an explicit $6$-$\pch$-critical caterpillar.

\begin{figure}[ht!]
\begin{center}
\begin{tikzpicture}[scale=0.4,style=thick]
\def\vr{5pt} 
\path (1,3) coordinate (v1); \path (5,3) coordinate (v2); \path (9,3) coordinate (v3); \path (13,3) coordinate (v4);
\path (17,3) coordinate (v5); \path (21,3) coordinate (v6); \path (25,3) coordinate (v7); \path (29,3) coordinate (v8);
\path (33,3) coordinate (v9);
\path (0,1) coordinate (w1); \path (1,1) coordinate (w2); \path (2,1) coordinate (w3); \path (4,1) coordinate (w4); \path (6,1) coordinate (w5);
\path (8,1) coordinate (w6); \path (10,1) coordinate (w7); \path (12,1) coordinate (w8); \path (13,1) coordinate (w9); \path (14,1) coordinate (w10);
\path (16,1) coordinate (w11); \path (18,1) coordinate (w12); \path (20,1) coordinate (w13); \path (21,1) coordinate (w14); \path (22,1) coordinate (w15);
\path (24,1) coordinate (w16); \path (26,1) coordinate (w17); \path (28,1) coordinate (w18); \path (30,1) coordinate (w19); \path (32,1) coordinate (w20);
\path (33,1) coordinate (w21); \path (34,1) coordinate (w22);

\draw (v1) -- (v2); \draw (v2) -- (v3); \draw (v3) -- (v4); \draw (v4) -- (v5); \draw (v5) -- (v6); \draw (v6) -- (v7);
\draw (v7) -- (v8); \draw (v8) -- (v9);

\draw (v1) -- (w1); \draw (v1) -- (w2); \draw (v1) -- (w3);
\draw (v2) -- (w4); \draw (v2) -- (w5);
\draw (v3) -- (w6); \draw (v3) -- (w7);
\draw (v4) -- (w8); \draw (v4) -- (w9); \draw (v4) -- (w10);
\draw (v5) -- (w11); \draw (v5) -- (w12);
\draw (v6) -- (w13); \draw (v6) -- (w14); \draw (v6) -- (w15);
\draw (v7) -- (w16); \draw (v7) -- (w17);
\draw (v8) -- (w18); \draw (v8) -- (w19);
\draw (v9) -- (w20); \draw (v9) -- (w21); \draw (v9) -- (w22);
%
\foreach \i in {1,...,9}
{  \draw (v\i)  [fill=white] circle (\vr); }
\foreach \i in {1,...,22}
{  \draw (w\i)  [fill=white] circle (\vr); }
\draw[anchor=south] (v1) node{$r$}; \draw[anchor=south] (v2) node{$s$}; \draw[anchor=south] (v3) node{$t$};
\draw[anchor=south] (v4) node{$u$}; \draw[anchor=south] (v5) node{$v$}; \draw[anchor=south] (v6) node{$w$};
\draw[anchor=south] (v7) node{$x$}; \draw[anchor=south] (v8) node{$y$}; \draw[anchor=south] (v9) node{$z$};

\end{tikzpicture}
\end{center}
\caption{The caterpillar $T$.} \label{fig:T}
\end{figure}
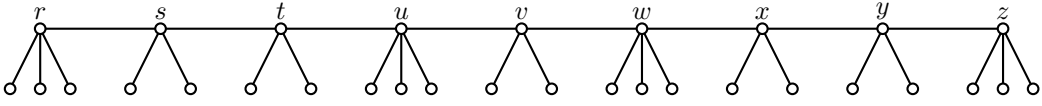

\begin{theorem}
\label{thm:caterpillar-chi=6}
The caterpillar $T$ from Fig.~\ref{fig:T} is $6$-$\pch$-critical.
\end{theorem}

To prove Theorem~\ref{thm:caterpillar-chi=6} we first demonstrate the following two lemmas.

\begin{lemma}
\label{lem:P9-no-(2,3,4,5)}
The path $P_9$ does not admit a $(2,3,4,5)$-coloring.
\end{lemma}

\proof
Let the vertices of $P_9$ be denoted with $r,s, \ldots, z$ in the natural order.
In a possible $(2,3,4,5)$-coloring of $P_9$, at most three vertices are colored $2$, at most three vertices are colored $3$, at most two vertices are colored $4$, and at
most two vertices are colored $5$. Suppose there exists a $(2,3,4,5)$-coloring $c:V(P_9) \rightarrow \{2,3,4,5\}$. If $|c^{-1}(3)| = 3$, then $c^{-1}(3) = \{r, v, z\}$.
However, this implies that $|c^{-1}(2)| = 2$ and $|c^{-1}(5)| = \{s, y\}$, which means at most one vertex from $\{t, u, w, x\}$ is colored $4$. Thus no such $c$
exists with $|c^{-1}(3)| = 3$.

Hence $|c^{-1}(2)| = 3$, $|c^{-1}(3)| = 2$, $|c^{-1}(4)| = 2 = |c^{-1}(5)|$. By symmetry, there are four cases to be considered.

\medskip\noindent
{\bf Case 1}. $c^{-1}(5) = \{r, x\}$. \\
The subpath $s, t, u, v, w$ has at most one vertex colored $4$ by $c$, which leads to a contradiction since the remaining four vertices cannot be colored using only
colors $2$ and $3$.

\medskip\noindent
{\bf Case 2}. $c^{-1}(5) = \{r, y\}$. \\
In this case $|c^{-1}(2)\cap \{s, \ldots, x\}|\le 2$, which implies $c(z) = 2$. This in turn forces $c^{-1}(4) = \{s, x\}$. But since
$|c^{-1}(3)\cap \{t,u, v, w\}|\le 1$, we get a contradiction.

\medskip\noindent
{\bf Case 3}. $c^{-1}(5) = \{r, z\}$. \\
This assumption together with $|c^{-1}(2)| = 3$ implies $c^{-1}(2) = \{s, v, y\}$, which forces $|c^{-1}(4)\cap \{t, u, w, x\}| = 2$, a contradiction.

\medskip\noindent
{\bf Case 4}. $c^{-1}(5) = \{s, y\}$. \\
Finally, $c^{-1}(5) = \{s, y\}$ implies that $c(r) = 4$ or $c(z) = 4$. Assume without loss of generality that $c(r) = 4$. If also $c(z) = 4$, then since
$|c^{-1}(2)\cap \{t, u, v, w, x\}| \le 2$, we get a contradiction. Hence $c(z) = 2$.  Since  $|c^{-1}(3)| = 2$, we must also have $c(t) = c(x) = 3$. But
then the coloring can clearly not be completed.
\qed

\begin{lemma}
\label{lem:pch(T)=6}
$\pch(T) = 6$.
\end{lemma}

\proof
Let the spine of $T$ be the set $S=\{r,s,t,u,v,w,x,y,z\}$, see Fig.~\ref{fig:T}. For any vertex $a$ in the spine $S$, if $a$ is adjacent to $k$ leaves, then for convenience we  denote this set of leaves by $\{a_1,\ldots,a_k\}$.

Coloring the vertices of the spine of $T$ with colors $2, 3, 4, 2, 5, 3, 2, 4,6$, and each of the leaves with $1$ demonstrates that $\pch(T)\le 6$. Thus it remains to prove that $\pch(T)\ge 6$.

Suppose for the sake of contradiction that
$\pch(T)\le 5$ and let $c$ be a $5$-packing coloring of $T$. Then, since $P_9$ does not admit a $(2,3,4,5)$-coloring, it is the case that some vertex of $S$
is colored 1.  The neighbors of this vertex of $S$ that is colored 1 receive  pairwise different colors from $\{2,3,4,5\}$. Consequently, the only
vertices of $S$ that can possibly be colored 1 are $r,s,t,v,x,y$ and $z$.  By symmetry we consider four cases, namely $c(r)=1$, $c(s)=1$, $c(t)=1$ and
$c(v)=1$.

\medskip\noindent
{\bf Case 1.} $c(r)=1$.  \\
 Since $c(\{r_1,r_2,r_3,s\})=\{2,3,4,5\}$ and the neighbors of $t$ are at distance at most 4 from the vertices in $\{r_1,r_2,r_3,s\}$,
it follows that $c(t) \not=1$. This implies that $c(t)=2$, that some vertex in $\{r_1,r_2,r_3\}$ is colored 2, and also that $c(u) \not=1$.  Considering the distance from $u$ to
the set $\{r_1,r_2,r_3,s\}$ we infer that $c(u)=3$ and also that $3 \in c(\{r_1,r_2,r_3\})$.  As a result either $c(v)=1$ or $c(v)=4$.  If $c(v)=1$,
then $c(\{v_1,v_2,w\})=\{2,4,5\}$, which implies that $5 \in c(\{r_1,r_2,r_3\})$ and thus $c(s)=4$.  However, this is a contradiction since the distance
between $s$ and any vertex in $\{v_1,v_2,w\}$ is 4.  Therefore, $c(v)=4$ and $4 \in c(\{r_1,r_2,r_3\})$.  It now follows that $c(s)=5$, $c(w)=2$ and $c(x)=1$.
This is a contradiction since the distance between $v$ and every vertex in $\{x_1,x_2,y\}$ is 3.

\medskip\noindent
{\bf Case 2.} $c(s)=1$.  \\
In this case  $c(\{r,s_1,s_2,t\})=\{2,3,4,5\}$.  This implies that $c(u)=2$, which implies $c(v)=3$ and $c(w)=4$. (Note that $c(v)\ne 1$ since every neighbor of
$v$ is within distance 5 of $\{r,s_1,s_2,t\}$.)
Now it follows that $c(t)=5$ and $c(x)\not=1$ since the distance between $v$ and any neighbor of $x$ is at most 3.  We now infer that
$c(x)=2$, which forces $c(y)=1$.  But this is a contradiction since the distance between $w$ and any vertex in $N(y)$ is at most 3.

\medskip\noindent
{\bf Case 3.}  $c(t)=1$.  \\
It follows that $c(\{s,t_1,t_2,u\})= \{2,3,4,5\}$.  By case 1 we may assume that $c(r)\not=1$, which implies that in fact $c(r)=2$, for otherwise the color $c(r)\in \{3,4,5\}$  could not be used on either of the vertices $s$, $t_1$, $t_2$, $u$.
If $c(u)=2$, then  because $c(\{s,t_1,t_2,u\})= \{2,3,4,5\}$ we must have $c(v)=1$.  But in this case one of $v_1,v_2$ or $w$ is colored 5, which contradicts the fact that some neighbor of vertex $t$ is also
colored 5.  If $c(u)\not=2$, then $c(v) \in \{1,2\}$.  Suppose first that $c(v)=2$. Then $c(w)\in \{1,3\}$, where the possibility $c(w) = 1$ is impossible because of the degree of $w$. So $c(w) = 3$. Then either $c(s) = 5$ and $c(u) = 4$, or $c(s) = 4$ and $c(u) = 5$. The first case forces $c(x) = 1$, but then one of the neighbors of $x$ should have been colored $4$, but this is not possible because $c(u) = 4$. And if $c(s) = 4$ and $c(u) = 5$, then $c(x)\ne 1$ because every neighbor of $x$ is of distance at most 4 from $u$, which is colored 5.  Hence $c(x) = 4$. But then $\{c(y), c(z)\} = \{1,2\}$  which readily leads to a contradiction since by symmetry and Cases 1 and 2 we may assume that $c(y)\ne 1$ and $c(z)\ne 1$.   On the other hand,
if $c(v)=1$,  then $\max\{c(v_1),c(v_2),c(w)\} \ge 4$, which again is in conflict with
the fact that the distance between any of $v_1,v_2$ or $w$ and any vertex in $\{s,t_1,t_2,u\}$ is at most 4.

\medskip\noindent
{\bf Case 4.} $c(v)=1$.  \\
We have shown in the first three cases that $1 \not\in c(\{r,s,t\})$ and by symmetry that  $1 \not\in c(\{x,y,z\})$.
Furthermore, $c(u)\not=1$ and $c(w)\not=1$ since they are neighbors of $v$. Hence, $c(\{r,s,t,u,v_1,v_2,w,x,y,z\})=\{2,3,4,5\}$.  Since  $c(\{u,w,t_1,t_2\})= \{2,3,4,5\}$, at most
three of the vertices in $\{r,s,t,u,v_1,v_2,w,x,y,z\}$ can be colored 2, at most three can be colored 3, at most three can be colored 4 and
at most one can be colored 5.  The only way that three of these ten vertices can be colored 4 is if $c(r)=4=c(z)$ and one of $v_1,v_2$ is colored 4.
But in this case it is easy to see that $|c^{-1}(\{2,3\})|\le 5$.  This gives the contradiction that $|c^{-1}(\{2,3,4,5\})| \le 9$.  However, if
$|c^{-1}(4)|\le 2$, then again we arrive at  $|c^{-1}(\{2,3,4,5\})| \le 9$.
\qed

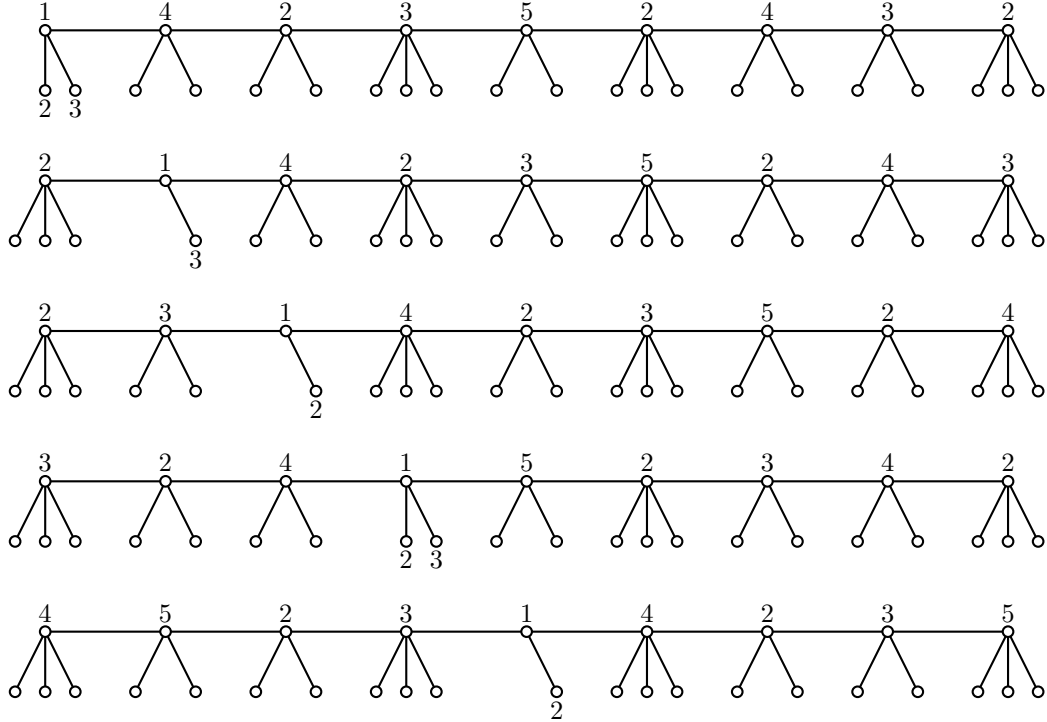
\begin{figure}[ht!]
\begin{center}
\begin{tikzpicture}[scale=0.4,style=thick]
\def\vr{5pt} 
\path (1,3) coordinate (va1); \path (5,3) coordinate (va2); \path (9,3) coordinate (va3); \path (13,3) coordinate (va4);
\path (17,3) coordinate (va5); \path (21,3) coordinate (va6); \path (25,3) coordinate (va7); \path (29,3) coordinate (va8);
\path (33,3) coordinate (va9);
\path (0,1) coordinate (wa1); \path (1,1) coordinate (wa2); \path (2,1) coordinate (wa3); \path (4,1) coordinate (wa4); \path (6,1) coordinate (wa5);
\path (8,1) coordinate (wa6); \path (10,1) coordinate (wa7); \path (12,1) coordinate (wa8); \path (13,1) coordinate (wa9); \path (14,1) coordinate (wa10);
\path (16,1) coordinate (wa11); \path (18,1) coordinate (wa12); \path (20,1) coordinate (wa13); \path (21,1) coordinate (wa14); \path (22,1) coordinate (wa15);
\path (24,1) coordinate (wa16); \path (26,1) coordinate (wa17); \path (28,1) coordinate (wa18); \path (30,1) coordinate (wa19); \path (32,1) coordinate (wa20);
\path (33,1) coordinate (wa21); \path (34,1) coordinate (wa22);

\draw (va1) -- (va2); \draw (va2) -- (va3); \draw (va3) -- (va4); \draw (va4) -- (va5); \draw (va5) -- (va6); \draw (va6) -- (va7);
\draw (va7) -- (va8); \draw (va8) -- (va9);

\draw (va1) -- (wa2); \draw (va1) -- (wa3);
\draw (va2) -- (wa4); \draw (va2) -- (wa5);
\draw (va3) -- (wa6); \draw (va3) -- (wa7);
\draw (va4) -- (wa8); \draw (va4) -- (wa9); \draw (va4) -- (wa10);
\draw (va5) -- (wa11); \draw (va5) -- (wa12);
\draw (va6) -- (wa13); \draw (va6) -- (wa14); \draw (va6) -- (wa15);
\draw (va7) -- (wa16); \draw (va7) -- (wa17);
\draw (va8) -- (wa18); \draw (va8) -- (wa19);
\draw (va9) -- (wa20); \draw (va9) -- (wa21); \draw (va9) -- (wa22);
%
\foreach \i in {1,...,9}
{  \draw (va\i)  [fill=white] circle (\vr); }
\foreach \i in {2,...,22}
{  \draw (wa\i)  [fill=white] circle (\vr); }
\draw[anchor=south] (va1) node{$1$}; \draw[anchor=south] (va2) node{$4$}; \draw[anchor=south] (va3) node{$2$};
\draw[anchor=south] (va4) node{$3$}; \draw[anchor=south] (va5) node{$5$}; \draw[anchor=south] (va6) node{$2$};
\draw[anchor=south] (va7) node{$4$}; \draw[anchor=south] (va8) node{$3$}; \draw[anchor=south] (va9) node{$2$};

\draw[anchor=north] (wa2) node{$2$}; \draw[anchor=north] (wa3) node{$3$};


\path (1,-2) coordinate (vb1); \path (5,-2) coordinate (vb2); \path (9,-2) coordinate (vb3); \path (13,-2) coordinate (vb4);
\path (17,-2) coordinate (vb5); \path (21,-2) coordinate (vb6); \path (25,-2) coordinate (vb7); \path (29,-2) coordinate (vb8);
\path (33,-2) coordinate (vb9);
\path (0,-4) coordinate (wb1); \path (1,-4) coordinate (wb2); \path (2,-4) coordinate (wb3); \path (4,-4) coordinate (wb4); \path (6,-4) coordinate (wb5);
\path (8,-4) coordinate (wb6); \path (10,-4) coordinate (wb7); \path (12,-4) coordinate (wb8); \path (13,-4) coordinate (wb9); \path (14,-4) coordinate (wb10);
\path (16,-4) coordinate (wb11); \path (18,-4) coordinate (wb12); \path (20,-4) coordinate (wb13); \path (21,-4) coordinate (wb14); \path (22,-4) coordinate (wb15);
\path (24,-4) coordinate (wb16); \path (26,-4) coordinate (wb17); \path (28,-4) coordinate (wb18); \path (30,-4) coordinate (wb19); \path (32,-4) coordinate (wb20);
\path (33,-4) coordinate (wb21); \path (34,-4) coordinate (wb22);

\draw (vb1) -- (vb2); \draw (vb2) -- (vb3); \draw (vb3) -- (vb4); \draw (vb4) -- (vb5); \draw (vb5) -- (vb6); \draw (vb6) -- (vb7);
\draw (vb7) -- (vb8); \draw (vb8) -- (vb9);

\draw (vb1) -- (wb1); \draw (vb1) -- (wb2); \draw (vb1) -- (wb3);
\draw (vb2) -- (wb5);
\draw (vb3) -- (wb6); \draw (vb3) -- (wb7);
\draw (vb4) -- (wb8); \draw (vb4) -- (wb9); \draw (vb4) -- (wb10);
\draw (vb5) -- (wb11); \draw (vb5) -- (wb12);
\draw (vb6) -- (wb13); \draw (vb6) -- (wb14); \draw (vb6) -- (wb15);
\draw (vb7) -- (wb16); \draw (vb7) -- (wb17);
\draw (vb8) -- (wb18); \draw (vb8) -- (wb19);
\draw (vb9) -- (wb20); \draw (vb9) -- (wb21); \draw (vb9) -- (wb22);
%
\foreach \i in {1,...,9}
{  \draw (vb\i)  [fill=white] circle (\vr); }
\foreach \i in {1,2,3,5,6,7,8,9,10,11,12,13,14,15,16,17,18,19,20,21,22}
{  \draw (wb\i)  [fill=white] circle (\vr); }
\draw[anchor=south] (vb1) node{$2$}; \draw[anchor=south] (vb2) node{$1$}; \draw[anchor=south] (vb3) node{$4$};
\draw[anchor=south] (vb4) node{$2$}; \draw[anchor=south] (vb5) node{$3$}; \draw[anchor=south] (vb6) node{$5$};
\draw[anchor=south] (vb7) node{$2$}; \draw[anchor=south] (vb8) node{$4$}; \draw[anchor=south] (vb9) node{$3$};

\draw[anchor=north] (wb5) node{$3$}; 


\path (1,-7) coordinate (vc1); \path (5,-7) coordinate (vc2); \path (9,-7) coordinate (vc3); \path (13,-7) coordinate (vc4);
\path (17,-7) coordinate (vc5); \path (21,-7) coordinate (vc6); \path (25,-7) coordinate (vc7); \path (29,-7) coordinate (vc8);
\path (33,-7) coordinate (vc9);
\path (0,-9) coordinate (wc1); \path (1,-9) coordinate (wc2); \path (2,-9) coordinate (wc3); \path (4,-9) coordinate (wc4); \path (6,-9) coordinate (wc5);
\path (8,-9) coordinate (wc6); \path (10,-9) coordinate (wc7); \path (12,-9) coordinate (wc8); \path (13,-9) coordinate (wc9); \path (14,-9) coordinate (wc10);
\path (16,-9) coordinate (wc11); \path (18,-9) coordinate (wc12); \path (20,-9) coordinate (wc13); \path (21,-9) coordinate (wc14); \path (22,-9) coordinate (wc15);
\path (24,-9) coordinate (wc16); \path (26,-9) coordinate (wc17); \path (28,-9) coordinate (wc18); \path (30,-9) coordinate (wc19); \path (32,-9) coordinate (wc20);
\path (33,-9) coordinate (wc21); \path (34,-9) coordinate (wc22);

\draw (vc1) -- (vc2); \draw (vc2) -- (vc3); \draw (vc3) -- (vc4); \draw (vc4) -- (vc5); \draw (vc5) -- (vc6); \draw (vc6) -- (vc7);
\draw (vc7) -- (vc8); \draw (vc8) -- (vc9);

\draw (vc1) -- (wc1); \draw (vc1) -- (wc2); \draw (vc1) -- (wc3);
\draw (vc2) -- (wc4); \draw (vc2) -- (wc5);
\draw (vc3) -- (wc7);
\draw (vc4) -- (wc8); \draw (vc4) -- (wc9); \draw (vc4) -- (wc10);
\draw (vc5) -- (wc11); \draw (vc5) -- (wc12);
\draw (vc6) -- (wc13); \draw (vc6) -- (wc14); \draw (vc6) -- (wc15);
\draw (vc7) -- (wc16); \draw (vc7) -- (wc17);
\draw (vc8) -- (wc18); \draw (vc8) -- (wc19);
\draw (vc9) -- (wc20); \draw (vc9) -- (wc21); \draw (vc9) -- (wc22);
%
\foreach \i in {1,...,9}
{  \draw (vc\i)  [fill=white] circle (\vr); }
\foreach \i in {1,2,3,4,5,7,8,9,10,11,12,13,14,15,16,17,18,19,20,21,22}
{  \draw (wc\i)  [fill=white] circle (\vr); }
\draw[anchor=south] (vc1) node{$2$}; \draw[anchor=south] (vc2) node{$3$}; \draw[anchor=south] (vc3) node{$1$};
\draw[anchor=south] (vc4) node{$4$}; \draw[anchor=south] (vc5) node{$2$}; \draw[anchor=south] (vc6) node{$3$};
\draw[anchor=south] (vc7) node{$5$}; \draw[anchor=south] (vc8) node{$2$}; \draw[anchor=south] (vc9) node{$4$};

\draw[anchor=north] (wc7) node{$2$}; 


\path (1,-12) coordinate (vd1); \path (5,-12) coordinate (vd2); \path (9,-12) coordinate (vd3); \path (13,-12) coordinate (vd4);
\path (17,-12) coordinate (vd5); \path (21,-12) coordinate (vd6); \path (25,-12) coordinate (vd7); \path (29,-12) coordinate (vd8);
\path (33,-12) coordinate (vd9);
\path (0,-14) coordinate (wd1); \path (1,-14) coordinate (wd2); \path (2,-14) coordinate (wd3); \path (4,-14) coordinate (wd4); \path (6,-14) coordinate (wd5);
\path (8,-14) coordinate (wd6); \path (10,-14) coordinate (wd7); \path (12,-14) coordinate (wd8); \path (13,-14) coordinate (wd9); \path (14,-14) coordinate (wd10);
\path (16,-14) coordinate (wd11); \path (18,-14) coordinate (wd12); \path (20,-14) coordinate (wd13); \path (21,-14) coordinate (wd14); \path (22,-14) coordinate (wd15);
\path (24,-14) coordinate (wd16); \path (26,-14) coordinate (wd17); \path (28,-14) coordinate (wd18); \path (30,-14) coordinate (wd19); \path (32,-14) coordinate (wd20);
\path (33,-14) coordinate (wd21); \path (34,-14) coordinate (wd22);

\draw (vd1) -- (vd2); \draw (vd2) -- (vd3); \draw (vd3) -- (vd4); \draw (vd4) -- (vd5); \draw (vd5) -- (vd6); \draw (vd6) -- (vd7);
\draw (vd7) -- (vd8); \draw (vd8) -- (vd9);

\draw (vd1) -- (wd1); \draw (vd1) -- (wd2); \draw (vd1) -- (wd3);
\draw (vd2) -- (wd4); \draw (vd2) -- (wd5);
\draw (vd3) -- (wd6); \draw (vd3) -- (wd7);
\draw (vd4) -- (wd9); \draw (vd4) -- (wd10);
\draw (vd5) -- (wd11); \draw (vd5) -- (wd12);
\draw (vd6) -- (wd13); \draw (vd6) -- (wd14); \draw (vd6) -- (wd15);
\draw (vd7) -- (wd16); \draw (vd7) -- (wd17);
\draw (vd8) -- (wd18); \draw (vd8) -- (wd19);
\draw (vd9) -- (wd20); \draw (vd9) -- (wd21); \draw (vd9) -- (wd22);
%
\foreach \i in {1,...,9}
{  \draw (vd\i)  [fill=white] circle (\vr); }
\foreach \i in {1,2,3,4,5,6,7,9,10,11,12,13,14,15,16,17,18,19,20,21,22}
{  \draw (wd\i)  [fill=white] circle (\vr); }
\draw[anchor=south] (vd1) node{$3$}; \draw[anchor=south] (vd2) node{$2$}; \draw[anchor=south] (vd3) node{$4$};
\draw[anchor=south] (vd4) node{$1$}; \draw[anchor=south] (vd5) node{$5$}; \draw[anchor=south] (vd6) node{$2$};
\draw[anchor=south] (vd7) node{$3$}; \draw[anchor=south] (vd8) node{$4$}; \draw[anchor=south] (vd9) node{$2$};

\draw[anchor=north] (wd9) node{$2$}; \draw[anchor=north] (wd10) node{$3$};


\path (1,-17) coordinate (ve1); \path (5,-17) coordinate (ve2); \path (9,-17) coordinate (ve3); \path (13,-17) coordinate (ve4);
\path (17,-17) coordinate (ve5); \path (21,-17) coordinate (ve6); \path (25,-17) coordinate (ve7); \path (29,-17) coordinate (ve8);
\path (33,-17) coordinate (ve9);
\path (0,-19) coordinate (we1); \path (1,-19) coordinate (we2); \path (2,-19) coordinate (we3); \path (4,-19) coordinate (we4); \path (6,-19) coordinate (we5);
\path (8,-19) coordinate (we6); \path (10,-19) coordinate (we7); \path (12,-19) coordinate (we8); \path (13,-19) coordinate (we9); \path (14,-19) coordinate (we10);
\path (16,-19) coordinate (we11); \path (18,-19) coordinate (we12); \path (20,-19) coordinate (we13); \path (21,-19) coordinate (we14); \path (22,-19) coordinate (we15);
\path (24,-19) coordinate (we16); \path (26,-19) coordinate (we17); \path (28,-19) coordinate (we18); \path (30,-19) coordinate (we19); \path (32,-19) coordinate (we20);
\path (33,-19) coordinate (we21); \path (34,-19) coordinate (we22);

\draw (ve1) -- (ve2); \draw (ve2) -- (ve3); \draw (ve3) -- (ve4); \draw (ve4) -- (ve5); \draw (ve5) -- (ve6); \draw (ve6) -- (ve7);
\draw (ve7) -- (ve8); \draw (ve8) -- (ve9);

\draw (ve1) -- (we1); \draw (ve1) -- (we2); \draw (ve1) -- (we3);
\draw (ve2) -- (we4); \draw (ve2) -- (we5);
\draw (ve3) -- (we6); \draw (ve3) -- (we7);
\draw (ve4) -- (we8); \draw (ve4) -- (we9); \draw (ve4) -- (we10);
\draw (ve5) -- (we12);
\draw (ve6) -- (we13); \draw (ve6) -- (we14); \draw (ve6) -- (we15);
\draw (ve7) -- (we16); \draw (ve7) -- (we17);
\draw (ve8) -- (we18); \draw (ve8) -- (we19);
\draw (ve9) -- (we20); \draw (ve9) -- (we21); \draw (ve9) -- (we22);
%
\foreach \i in {1,...,9}
{  \draw (ve\i)  [fill=white] circle (\vr); }
\foreach \i in {1,2,3,4,5,6,7,8,9,10,12,13,14,15,16,17,18,19,20,21,22}
{  \draw (we\i)  [fill=white] circle (\vr); }
\draw[anchor=south] (ve1) node{$4$}; \draw[anchor=south] (ve2) node{$5$}; \draw[anchor=south] (ve3) node{$2$};
\draw[anchor=south] (ve4) node{$3$}; \draw[anchor=south] (ve5) node{$1$}; \draw[anchor=south] (ve6) node{$4$};
\draw[anchor=south] (ve7) node{$2$}; \draw[anchor=south] (ve8) node{$3$}; \draw[anchor=south] (ve9) node{$5$};

\draw[anchor=north] (we12) node{$2$}; 

\end{tikzpicture}
\end{center}
\caption{5-packing colorings of $T-r_1$,  of $T-s_1$, of $T-t_1$, of $T-u_1$, and of $T-v_1$.}
\label{fig:Tr1}
\end{figure}

To complete the proof of Theorem~\ref{thm:caterpillar-chi=6} we must show that the removal of any vertex of $T$ leaves a graph with packing chromatic number at
most 5.  By using the $(2,3,4,5,6)$-coloring of $P_9$ given above it follows that any caterpillar whose spine has order at most 8
has packing chromatic number at most 5.  Consequently, we only need to show that there is a 5-packing coloring
of the caterpillar that remains when one of the leaves of $T$ is deleted.  Using the symmetry of $T$, such  packing colorings are shown in a sequence of five figures.
See Fig.~\ref{fig:Tr1}. In each case the leaves that are not labeled are colored 1.  Theorem~\ref{thm:caterpillar-chi=6} is thus proved.

\section{Critical Cartesian products}
\label{sec:Cartesian products}

The Cartesian product $G\cp H$ of graphs $G$ and $H$ has $V(G)\times V(H)$ as the vertex set, vertices $(g,h)$ and $(g',h')$ being adjacent if either $gg'\in E(G)$ and $h=h'$, or $g=g'$ and $hh'\in E(H)$, see~\cite{imrich-2008}.

Recall that a graph $G$ is vertex-critical for the chromatic number, if $\chi(G-u) = \chi(G) -1$ holds for every $u\in V(G)$. If the factors of a Cartesian product are both non-trivial, then $G\cp H$ is not vertex-critical for the chromatic number. Indeed, recall that $\chi(G\cp H) = \max\{\chi(G), \chi(H)\}$, see~\cite[Theorem 8.1]{imrich-2008}). Hence, since a vertex-deleted subgraph of $G\cp H$ contains a subgraph isomorphic to $G$ and a subgraph isomorphic to $H$, the vertex-deleted subgraph has the same chromatic number as $G\cp H$. In this section we show that on the other hand the Cartesian product is a very rich source of $\pch$-critical graphs. For additional investigation of $\pch$ on Cartesian products see~\cite{bkr-2007} and \cite{jacobs-2013}.

We begin with the following sporadic, but very illustrative example. Let $G=K_{1,3}\cp P_3$ with the vertices labeled as in Fig.~\ref{f:critical-CP}.  We claim that $G$ is $6$-$\pch$-critical even
though neither $K_{1,3}$ nor $P_3$ is critical.  First note that the independence number of $G$ is 7, and that $G$ has a unique independent
set of cardinality 7, namely $I=\{u_1,u_2,u_3,v,w_1,w_2,w_3\}$.  Furthermore, $\diam(G)=4$ and so in any packing coloring of $G$ any color
4 or more can be assigned at most once.

\begin{figure}[ht!]
\tikzstyle{every node}=[circle, draw, fill=black!0, inner sep=0pt,minimum width=.2cm]
\begin{center}
\begin{tikzpicture}[thick,scale=.7]
  \draw(0,0) { 
    +(0.00,4.00) -- +(1.25,4.00)
    +(1.25,4.00) -- +(2.50,4.00)
    +(1.25,4.00) -- +(3.75,5.00)
    +(3.75,3.00) -- +(1.25,2.00)
    +(1.25,2.00) -- +(2.50,2.00)
    +(1.25,2.00) -- +(0.00,2.00)
    +(0.00,0.00) -- +(1.25,0.00)
    +(1.25,0.00) -- +(3.75,1.00)
    +(1.25,0.00) -- +(2.50,0.00)
    +(3.75,5.00) -- +(3.75,3.00)
    +(3.75,3.00) -- +(3.75,1.00)
    +(2.50,4.00) -- +(2.50,2.00)
    +(2.50,2.00) -- +(2.50,0.00)
    +(1.25,4.00) -- +(1.25,2.00)
    +(1.25,2.00) -- +(1.25,0.00)
    +(0.00,4.00) -- +(0.00,2.00)
    +(0.00,2.00) -- +(0.00,0.00)
    +(0.00,4.00) node{}
    +(1.25,4.00) node{}
    +(2.50,4.00) node{}
    +(3.75,5.00) node{}
    +(0.00,2.00) node{}
    +(1.25,2.00) node{}
    +(2.50,2.00) node{}
    +(3.75,3.00) node{}
    +(0.00,0.00) node{}
    +(1.25,0.00) node{}
    +(2.50,0.00) node{}
    +(3.75,1.00) node{}
    +(4.3,1) node[rectangle, draw=white!0, fill=white!100]{$w_3$}
    +(4.3,3) node[rectangle, draw=white!0, fill=white!100]{$v_3$}
    +(4.3,5) node[rectangle, draw=white!0, fill=white!100]{$u_3$}
    +(3.05,0) node[rectangle, draw=white!0, fill=white!100]{$w_2$}
    +(3.05,2) node[rectangle, draw=white!0, fill=white!100]{$v_2$}
    +(3.05,4) node[rectangle, draw=white!0, fill=white!100]{$u_2$}
    +(-.55,0) node[rectangle, draw=white!0, fill=white!100]{$w_1$}
    +(-.55,2) node[rectangle, draw=white!0, fill=white!100]{$v_1$}
    +(-.55,4) node[rectangle, draw=white!0, fill=white!100]{$u_1$}
    +(.9,.3) node[rectangle, draw=white!0, fill=white!100]{$w$}
    +(.9,2.3) node[rectangle, draw=white!0, fill=white!100]{$v$}
    +(.9,4.3) node[rectangle, draw=white!0, fill=white!100]{$u$}
  };
\end{tikzpicture}
\end{center}
\caption{$K_{1,3}\cp P_3$.} \label{f:critical-CP}
\end{figure}

The packing coloring in which each vertex of $I$ is colored $1$ and the five neighbors of $v$ are assigned distinct colors
from the set $\{2,3,4,5,6\}$ shows that $\pch(G)\le 6$.  To prove that
$\pch(G) \ge 6$ we assume that $c: V(G) \longrightarrow [n]$ is a packing coloring with $n$ as small as possible.
If $c(v)=1$, then all of the vertices in $N(v)$ have distinct values under $c$, and this implies that $n\ge 6$.  If $c(v)=2$, then
since the eccentricity of $v$ is 2, it follows that $c^{-1}(2)=\{v\}$.  This implies that $|c^{-1}(1)|\le 6$.  In this case
if $|c^{-1}(3)| \ge 2$, then $3 \in c(I-\{v\})$.  Consequently, $|c^{-1}(\{1,3\})| \le 7$ and hence in this case
$n\ge 7$.  Similarly,  if $c(v)=3$, then $c^{-1}(3)=\{v\}$ and $|c^{-1}(\{1,2\})| \le 7$, which again means $n\ge 7$.  Finally,
if $c(v) \in \{4,5\}$, then one can easily show that $|c^{-1}(\{1,2,3\})| \le 9$, which leads to $n \ge 3+(12-9)=6$.  (One such
partial packing coloring $c$ with $|c^{-1}(\{1,2,3\})|= 9$ is shown in Table~\ref{tab:critical}.)

\begin{table}
\centering
\begin{tabular}{|c|cccccccccccc|} %
\hline
$x$&$u$&$u_1$&$u_2$&$u_3$&$v$&$v_1$&$v_2$&$v_3$&$w$&$w_1$&$w_2$&$w_3$  \\ \hline
$c(x)$ & &3&1&2&4&1&2&1&1&2& &3 \\ \hline
$f_u(x)$& &3&3&2&5&1&1&1&1&2&4&3 \\ \hline
$f_{u_3}(x)$ & 1&3&2& & 5&1&1&1&1&2&3&4 \\ \hline
$f_{v_3}(x)$ & 2&1&1&1&1&3&5& &4&1&1&1 \\ \hline
$f_{v}(x)$ &1&3&5&2& & 1&1&1&1&2&4&3\\ \hline
\end{tabular}
\caption{Some partial packing colorings of $K_{1,3}\cp P_3$.}
\label{tab:critical}
\end{table}

Finally to show that $G$ is critical, we observe that by symmetry we need to show that $\pch(G-t)<6$ for each vertex $t \in \{u,v,u_3,v_3\}$.
In the table there is a packing coloring, $f_t: V(G-t) \longrightarrow [5]$ for each such vertex $t$.  This proves that
$G$ is $6$-critical.

Note that $P_3$ and $K_{1,3}$ are not critical. Hence the above example shows that even the Cartesian product of non-critical graphs can be critical.

To establish a greater variety of $\pch$-critical Cartesian products, we recall an earlier result and prove a new lemma that might be of independent interest.

\begin{theorem} {\rm \cite[Theorem 1]{bkr-2007}}
\label{thm:lower-cart}
If $G$ and $H$ are connected graphs on at least two vertices, then
$$\chi_{\rho}(G\cp H) \geq (\chi_{\rho}(G)+1)n(H) - {\rm diam}(G\cp H)(n(H)- 1) -1\,.$$
\end{theorem}

\begin{lemma}
\label{lem:vt}
If $G$ is a vertex-transitive graph with $\diam(G) \le \pch(G)$, then $G$ is $\pch$-critical.
\end{lemma}

\proof
Let $x$ be any vertex of $G$.  Since $G$ is vertex-transitive, there is a $\pch(G)$-packing coloring $c$ of $G$ such that $c(x)=\pch(G)$.  By Observation~\ref{obs:largecolors}, $x$ is the only vertex assigned the color $\pch(G)$ by $c$.  The restriction of $c$ to $V(G-x)$ is a packing coloring, which shows that $\pch(G-x)<\pch(G)$.
\qed

Consider the hypercubes $Q_n$, $n\ge 1$. Indeed, $\pch(Q_n)$ was determined exactly for $n\le 5$ in~\cite{goddard-2008} and for $6\le n\le 8$ in~\cite{torres-2015}. Moreover, $\pch(Q_n)$ asymptotically grows as $\left(\frac{1}{2} - O(\frac{1}{n})\right)2^n$. Since $\diam(Q_n) = n$, we conclude from Lemma~\ref{lem:vt} that $Q_n$ is critical for every $n\ge 1$.

The announced result now reads as follows.

\begin{theorem}
\label{thm:critical-products}
If $G$ and $H$ are connected, vertex-transitive graphs on at least two vertices and $\diam(G) + \diam(H) \le \pch(G)$, then $G\cp H$ is $\pch$-critical.
\end{theorem}

\proof
Since $\diam(G) + \diam(H) \le \pch(G)$ we have $\pch(G) - \diam(G) - \diam(H) + 1\ge 1$, from which we get that
\begin{equation}
\label{eq:at-least-1}
n(H)(\pch(G) - \diam(G) - \diam(H) + 1) \ge 1\,.
\end{equation}
Using Theorem~\ref{thm:lower-cart} and the well-known fact that $\diam(G\cp H) = \diam (G) + \diam(H)$ (see~\cite[p.102]{imrich-2008}), we can estimate as follows:
\begin{eqnarray*}
\pch(G \cp H) & \ge & (\chi_{\rho}(G)+1)n(H) - {\rm diam}(G\cp H)(n(H)- 1) - 1 \\
& = & n(H)(\pch(G) - \diam(G) - \diam(H) + 1) + \diam(G) + \diam(H) - 1 \\
& \ge & \diam(G) + \diam(H)\,,
\end{eqnarray*}
where the last inequality follows by~\eqref{eq:at-least-1}.

It is also well-known (see~\cite[Proposition 6.16]{hammack-2011}) that a Cartesian product of connected graphs is vertex-transitive if and only if the factors are such.
Hence, $G\cp H$ is a vertex-transitive graph. Moreover, by the above,
$$\diam(G\cp H) = \diam(G) + \diam(H) \le \pch(G \cp H)\,.$$
The result now follows by Lemma~\ref{lem:vt}.
\qed

Let $k\ge 3$ and let $H$ be a vertex-transitive graph with $\diam(H) \le k-1$. Then Theorem~\ref{thm:critical-products} implies that $K_k\cp H$ is $\pch$-critical. This fact in
particular implies that every vertex-transitive graph is an induced (actually convex) subgraph of a $\pch$-critical graph.

For a particular example consider $H = C_{4\ell}$, where $2\ell\le k-1$. Then $K_k\cp C_{4\ell}$ is $\pch$-critical graph, where one factor, namely $C_{4\ell}$, is not $\pch$-critical.

\begin{figure}[ht!]
\begin{center}
\begin{tikzpicture}[scale=0.75,style=thick]
\def\vr{2pt} 
\path (-9,0) coordinate (v1); \path (-9,1) coordinate (w1); \path (-8,0) coordinate (v2); \path (-8,1) coordinate (w2);
\path (-7,0) coordinate (v3); \path (-7,1) coordinate (w3); \path (-6,0) coordinate (v4); \path (-6,1) coordinate (w4);
\path (-5,0) coordinate (v5); \path (-5,1) coordinate (w5); \path (-4,0) coordinate (v6); \path (-4,1) coordinate (w6);
\path (-3,0) coordinate (v7); \path (-3,1) coordinate (w7); \path (-2,0) coordinate (v8); \path (-2,1) coordinate (w8);
\path (-1,0) coordinate (v9); \path (-1,1) coordinate (w9); \path (0,0) coordinate (v10); \path (0,1) coordinate (w10);
\path (1,0) coordinate (v11); \path (1,1) coordinate (w11); \path (2,0) coordinate (v12); \path (2,1) coordinate (w12);
\path (3,0) coordinate (v13); \path (3,1) coordinate (w13); \path (4,0) coordinate (v14); \path (4,1) coordinate (w14);
\path (5,0) coordinate (v15); \path (5,1) coordinate (w15); \path (6,0) coordinate (v16); \path (6,1) coordinate (w16);
\path (7,0) coordinate (v17); \path (7,1) coordinate (w17); \path (8,0) coordinate (v18); \path (8,1) coordinate (w18);
%
\draw (v1) -- (v2); \draw (w1) -- (w2);
\foreach \i in {1,...,18}
{  \draw (v\i) -- (w\i); }
\draw (v1) -- (v2); \draw (v2) -- (v3); \draw (v3) -- (v4); \draw (v4) -- (v5); \draw (v5) -- (v6); \draw (v6) -- (v7);
\draw (v7) -- (v8); \draw (v8) -- (v9); \draw (v9) -- (v10); \draw (v10) -- (v11); \draw (v11) -- (v12); \draw (v12) -- (v13);
\draw (v13) -- (v14); \draw (v14) -- (v15); \draw (v15) -- (v16); \draw (v16) -- (v17); \draw (v17) -- (v18);
\draw (w1) -- (w2); \draw (w2) -- (w3); \draw (w3) -- (w4); \draw (w4) -- (w5); \draw (w5) -- (w6); \draw (w6) -- (w7);
\draw (w7) -- (w8); \draw (w8) -- (w9); \draw (w9) -- (w10); \draw (w10) -- (w11); \draw (w11) -- (w12); \draw (w12) -- (w13);
\draw (w13) -- (w14); \draw (w14) -- (w15); \draw (w15) -- (w16); \draw (w16) -- (w17); \draw (w17) -- (w18);
\draw (v1) .. controls (-8,-2)and (7,-2) .. (v18);
\draw (w1) .. controls (-8,3)and (7,3) .. (w18);
%
\foreach \i in {1,...,18}
{  \draw (v\i)  [fill=white] circle (\vr); }
\foreach \i in {1,...,18}
{  \draw (w\i)  [fill=white] circle (\vr); }

\foreach \j in {1,3,5,7,9,11,13,15,17}
{   \draw[anchor = south] (w\j) node {$1$};  }
\foreach \j in {4,10,16}
{   \draw[anchor = south] (w\j) node {$2$};  }
\foreach \j in {2,8,14}
{   \draw[anchor = south] (w\j) node {$3$};  }
\foreach \j in {6,12,18}
{   \draw[anchor = south] (w\j) node {$5$};  }

\foreach \j in {2,4,6,8,10,12,14,16,18}
{   \draw[anchor = north] (v\j) node {$1$};  }
\foreach \j in {1,7,13}
{   \draw[anchor = north] (v\j) node {$2$};  }
\foreach \j in {5,11,17}
{   \draw[anchor = north] (v\j) node {$3$};  }
\foreach \j in {3,9,15}
{   \draw[anchor = north] (v\j) node {$4$};  }

\end{tikzpicture}
\end{center}
\caption{$\pch$-coloring of $C_{18} \cp P_2$}
\label{fig:ncproduct}
\end{figure}

Although many pairs of vertex-transitive graphs that do not satisfy the diameter condition in Theorem~\ref{thm:critical-products} still
produce a critical Cartesian product, it is not always the case.  See Fig.~\ref{fig:ncproduct} where a coloring of the product
$C_{18} \cp P_2$ is shown.  From \cite[Proposition 6.1]{goddard-2008} we know that the $2 \times k$ grid has packing chromatic number 5 for every $k \ge 6$.
We infer that $\pch(C_{18} \cp P_2) \ge 5$ since $P_{18} \cp P_2$ is a subgraph of $C_{18} \cp P_2$.  It is easy to check that the coloring in Fig.~\ref{fig:ncproduct}
is indeed a packing coloring, which in turn shows that $\pch(C_{18} \cp P_2) = 5$.  For any vertex $x$ of $C_{18} \cp P_2$, the $2 \times 17$ grid is a subgraph of $(C_{18} \cp P_2)-x$.  Since $C_{18} \cp P_2$ is vertex-transitive, this shows that $C_{18} \cp P_2$ is not $\pch$-critical.

\section{Concluding remarks}
\label{sec:conclude}

In this paper we have investigated $\pch$-critical graphs. These graphs are vertex critical for the packing chromatic number. An equally legal concept would be edge critical graphs, that is, graphs $G$ for which $\pch(G-e) < \pch(G)$ holds for every edge $e$ of $G$. Their investigation seems quite different from the vertex-criticality studied in this paper though.

It would be interesting to classify vertex-transitive, $\pch$-critical graphs. Examples include cycles (only those whose order is not congruent to $0$ modulo $4$), complete graphs, complete multipartite graphs with parts of equal size, hypercubes, the Petersen graph $P$, as well as the Cartesian product of those that satisfy the additional condition of Theorem~\ref{thm:critical-products}.
A sporadic example from the latter family of graphs is $P\cp P$.

We did not give a complete list of $4$-$\pch$-critical graphs. A different approach to this problem would be to use the structure of graphs $G$ with $\pch(G) = 3$ as described in~\cite{goddard-2008}. In particular, for the $2$-connected case we have the following characterization.

\begin{proposition} {\rm \cite[Proposition 3.2]{goddard-2008}} \label{prop:3critical2conn}
Let $G$ be a $2$-connected graph.  Then $\pch(G)=3$ if and only if $G$ is either the subdivision of a bipartite multigraph
or the join of $K_2$ and an independent set.
\end{proposition}

Based on this proposition we say that a graph $G$ is {\em join reducible} if there exists at least one vertex $x$ of $G$ such that $G-x$ is the join of $K_2$ and an independent set. Then we can prove the following result, the proof of which is a tedious case analysis and hence is omitted.

\begin{proposition}
If $G$ is a $4$-$\pch$-critical graph, then $G$ is join reducible if and only if $G$ is one of the five graphs in Fig.~\ref{fig:4critical}.
\end{proposition}

\begin{figure}[ht!]
\tikzstyle{every node}=[circle, draw, fill=black!0, inner sep=0pt,minimum width=.2cm]
\begin{center}
\begin{tikzpicture}[thick,scale=1.0]
  \draw(0,0) { 
    +(1.50,4.50) -- +(0.7,4)
    +(0.70,4.00) -- +(1.00,3.00)
    +(1.00,3.00) -- +(2.00,3.00)
    +(0.70,4.00) -- +(2.00,3.00)
    +(1.50,4.50) -- +(1.00,3.00)
    +(5.30,4.00) -- +(4.00,3.00)
    +(4.00,3.00) -- +(5.00,3.00)
    +(5.00,3.00) -- +(3.70,4.00)
    +(3.70,4.00) -- +(4.00,3.00)
    +(5.30,4.00) -- +(5.00,3.00)
    +(4.50,4.50) -- +(3.70,4.00)
    +(-0.30,1.00) -- +(0.00,0.00)
    +(0.00,0.00) -- +(1.00,0.00)
    +(1.00,0.00) -- +(1.30,1.00)
    +(1.30,1.00) -- +(0.00,0.00)
    +(-0.30,1.00) -- +(1.00,0.00)
    +(0.50,1.50) -- +(-0.30,1.00)
    +(0.50,1.50) -- +(1.30,1.00)
    +(2.50,0.00) -- +(3.50,0.00)
    +(3.50,0.00) -- +(2.20,1.00)
    +(2.20,1.00) -- +(2.50,0.00)
    +(2.50,0.00) -- +(3.80,1.00)
    +(3.80,1.00) -- +(3.50,0.00)
    +(3.00,1.50) -- +(2.20,1.00)
    +(3.00,1.50) -- +(2.50,0.00)
    +(1.50,4.50) -- +(2.00,3.00)
    +(5.00,0.00) -- +(6.00,0.00)
    +(6.00,0.00) -- +(4.70,1.00)
    +(4.70,1.00) -- +(5.00,0.00)
    +(5.00,0.00) -- +(6.30,1.00)
    +(6.30,1.00) -- +(6.00,0.00)
    +(5.50,1.50) -- +(4.70,1.00)
    +(6.30,1.00) -- +(5.50,1.50)
    +(5.50,1.50) -- +(5.00,0.00)
    +(0.7,4.00) node{}
    +(1.00,3.00) node{}
    +(2.00,3.00) node{}
    +(1.50,4.50) node{}
    +(4.00,3.00) node{}
    +(5.00,3.00) node{}
    +(3.7,4.00) node{}
    +(4.5,4.50) node{}
    +(5.30,4.00) node{}
    +(0.50,1.50) node{}
    +(-0.30,1.00) node{}
    +(1.30,1.00) node{}
    +(0.00,0.00) node{}
    +(1.00,0.00) node{}
    +(3.00,1.50) node{}
    +(2.20,1.00) node{}
    +(3.80,1.00) node{}
    +(2.50,0.00) node{}
    +(3.50,0.00) node{}
    +(4.70,1.00) node{}
    +(6.30,1.00) node{}
    +(5.00,0.00) node{}
    +(6.00,0.00) node{}
    +(5.50,1.50) node{}
    +(1.2,4.5) node[rectangle, draw=white!0, fill=white!100]{$x$}
    +(4.2,4.5) node[rectangle, draw=white!0, fill=white!100]{$x$}
    +(0.2,1.5) node[rectangle, draw=white!0, fill=white!100]{$x$}
    +(2.7,1.5) node[rectangle, draw=white!0, fill=white!100]{$x$}
    +(5.2,1.5) node[rectangle, draw=white!0, fill=white!100]{$x$}
  };
\end{tikzpicture}
\end{center}
\caption{Some $4$-$\pch$-critical graphs.} \label{fig:4critical}
\end{figure}

\nocite{*}
\bibliographystyle{abbrvnat}
\bibliography{biblio-packing-final}

\begin{thebibliography}{20}
\providecommand{\natexlab}[1]{#1}
\providecommand{\url}[1]{\texttt{#1}}
\expandafter\ifx\csname urlstyle\endcsname\relax
  \providecommand{\doi}[1]{doi: #1}\else
  \providecommand{\doi}{doi: \begingroup \urlstyle{rm}\Url}\fi

\bibitem[Balogh et~al.(2018)Balogh, Kostochka, and Liu]{balogh-2018}
J.~Balogh, A.~Kostochka, and X.~Liu.
\newblock Packing chromatic number of cubic graphs.
\newblock \emph{Discrete Math.}, 341\penalty0 (2):\penalty0 474--483, 2018.
\newblock \doi{10.1016/j.disc.2017.09.014}.

\bibitem[Bre\v{s}ar and Ferme(2018)]{bresar-2018a}
B.~Bre\v{s}ar and J.~Ferme.
\newblock An infinite family of subcubic graphs with unbounded packing
  chromatic number.
\newblock \emph{Discrete Math.}, 341\penalty0 (8):\penalty0 2337--2342, 2018.
\newblock \doi{10.1016/j.disc.2018.05.004}.

\bibitem[Bre\v{s}ar et~al.()Bre\v{s}ar, Gastineau, and Togni]{bresar-2018+}
B.~Bre\v{s}ar, N.~Gastineau, and O.~Togni.
\newblock Packing colorings of subcubic outerplanar graphs.
\newblock arXiv:1809.05552 [math.CO], 14 Sep 2018.

\bibitem[Bre\v{s}ar et~al.(2007)Bre\v{s}ar, Klav\v{z}ar, and Rall]{bkr-2007}
B.~Bre\v{s}ar, S.~Klav\v{z}ar, and D.~F. Rall.
\newblock On the packing chromatic number of {C}artesian products, hexagonal
  lattice, and trees.
\newblock \emph{Discrete Appl. Math.}, 155\penalty0 (17):\penalty0 2303--2311,
  2007.
\newblock \doi{10.1016/j.dam.2007.06.008}.

\bibitem[Bre\v{s}ar et~al.(2017)Bre\v{s}ar, Klav\v{z}ar, Rall, and
  Wash]{bkrw-2017a}
B.~Bre\v{s}ar, S.~Klav\v{z}ar, D.~F. Rall, and K.~Wash.
\newblock Packing chromatic number under local changes in a graph.
\newblock \emph{Discrete Math.}, 340\penalty0 (5):\penalty0 1110--1115, 2017.
\newblock \doi{10.1016/j.disc.2016.09.030}.

\bibitem[Bre\v{s}ar et~al.(2018)Bre\v{s}ar, Klav\v{z}ar, Rall, and
  Wash]{bresar-2018b}
B.~Bre\v{s}ar, S.~Klav\v{z}ar, D.~F. Rall, and K.~Wash.
\newblock Packing chromatic number versus chromatic and clique number.
\newblock \emph{Aequationes Math.}, 92\penalty0 (3):\penalty0 497--513, 2018.
\newblock \doi{10.1007/s00010-017-0520-9}.

\bibitem[Fiala and Golovach(2010)]{fiala-2010}
J.~Fiala and P.~A. Golovach.
\newblock Complexity of the packing coloring problem for trees.
\newblock \emph{Discrete Appl. Math.}, 158\penalty0 (7):\penalty0 771--778,
  2010.
\newblock \doi{10.1016/j.dam.2008.09.001}.

\bibitem[Gastineau and Togni(2016)]{gt-2016}
N.~Gastineau and O.~Togni.
\newblock {$S$}-packing colorings of cubic graphs.
\newblock \emph{Discrete Math.}, 339\penalty0 (10):\penalty0 2461--2470, 2016.
\newblock \doi{10.1016/j.disc.2016.04.017}.

\bibitem[Gastineau et~al.(2018)Gastineau, Holub, and Togni]{gastineau-2018}
N.~Gastineau, P.~Holub, and O.~Togni.
\newblock On the packing chromatic number of subcubic outerplanar graphs.
\newblock \emph{Discrete Appl. Math.}, 2018.
\newblock \doi{10.1016/j.dam.2018.07.034}.

\bibitem[Goddard and Xu(2012)]{goddard-2012}
W.~Goddard and H.~Xu.
\newblock The {$S$}-packing chromatic number of a graph.
\newblock \emph{Discuss. Math. Graph Theory}, 32\penalty0 (4):\penalty0
  795--806, 2012.
\newblock \doi{10.7151/dmgt.1642}.

\bibitem[Goddard et~al.(2008)Goddard, Hedetniemi, Hedetniemi, Harris, and
  Rall]{goddard-2008}
W.~Goddard, S.~M. Hedetniemi, S.~T. Hedetniemi, J.~M. Harris, and D.~F. Rall.
\newblock Broadcast chromatic numbers of graphs.
\newblock \emph{Ars Combin.}, 86:\penalty0 33--49, 2008.

\bibitem[Hammack et~al.(2011)Hammack, Imrich, and Klav\v{z}ar]{hammack-2011}
R.~Hammack, W.~Imrich, and S.~Klav\v{z}ar.
\newblock \emph{Handbook of product graphs}.
\newblock Discrete Mathematics and its Applications (Boca Raton). CRC Press,
  Boca Raton, FL, second edition, 2011.
\newblock With a foreword by Peter Winkler.

\bibitem[Imrich et~al.(2008)Imrich, Klav\v{z}ar, and Rall]{imrich-2008}
W.~Imrich, S.~Klav\v{z}ar, and D.~F. Rall.
\newblock \emph{Topics in graph theory}.
\newblock A K Peters, Ltd., Wellesley, MA, 2008.
\newblock Graphs and their Cartesian product.

\bibitem[Jacobs et~al.(2013)Jacobs, Jonck, and Joubert]{jacobs-2013}
Y.~Jacobs, E.~Jonck, and E.~J. Joubert.
\newblock A lower bound for the packing chromatic number of the {C}artesian
  product of cycles.
\newblock \emph{Cent. Eur. J. Math.}, 11\penalty0 (7):\penalty0 1344--1357,
  2013.
\newblock \doi{10.2478/s11533-013-0237-5}.

\bibitem[Jensen(2002)]{jensen-2002}
T.~R. Jensen.
\newblock Dense critical and vertex-critical graphs.
\newblock \emph{Discrete Math.}, 258\penalty0 (1-3):\penalty0 63--84, 2002.
\newblock \doi{10.1016/S0012-365X(02)00262-5}.

\bibitem[Kor\v{z}e and Vesel(2019)]{korze-2019}
D.~Kor\v{z}e and A.~Vesel.
\newblock Packing coloring of generalized {S}ierpi\'{n}ski graphs.
\newblock \emph{Discrete Math. Theor. Comput. Sci.}, 21\penalty0 (3):\penalty0
  Paper No. 7, 18, 2019.

\bibitem[La\"{i}che and Sopena(2018)]{laiche-2018}
D.~La\"{i}che and E.~Sopena.
\newblock Packing colouring of some classes of cubic graphs.
\newblock \emph{Australas. J. Combin.}, 72:\penalty0 376--404, 2018.

\bibitem[Martin et~al.(2017)Martin, Raimondi, Chen, and Martin]{barnaby-2017}
B.~Martin, F.~Raimondi, T.~Chen, and J.~Martin.
\newblock The packing chromatic number of the infinite square lattice is
  between 13 and 15.
\newblock \emph{Discrete Appl. Math.}, 225:\penalty0 136--142, 2017.
\newblock \doi{10.1016/j.dam.2017.03.013}.

\bibitem[Sloper(2004)]{sloper-2004}
C.~Sloper.
\newblock An eccentric coloring of trees.
\newblock \emph{Australas. J. Combin.}, 29:\penalty0 309--321, 2004.

\bibitem[Torres and Valencia-Pabon(2015)]{torres-2015}
P.~Torres and M.~Valencia-Pabon.
\newblock The packing chromatic number of hypercubes.
\newblock \emph{Discrete Appl. Math.}, 190/191:\penalty0 127--140, 2015.
\newblock \doi{10.1016/j.dam.2015.04.006}.

\end{thebibliography}
\label{sec:biblio}

\end{document}